\documentclass[12pt]{article}
 
\usepackage{bbm, amsmath, amsfonts, amscd, latexsym, amsthm, amssymb, graphicx}
\usepackage[nottoc,notlot,notlof]{tocbibind}
\newtheorem{theor}{\hspace{1cm}{\sc Theorem}}[section]
\newtheorem{utver}[theor]{\hspace{1cm}{\sc Proposition}}
\newtheorem{sledst}[theor]{\hspace{1cm}{\sc Corollary}}
\newtheorem{lemma}[theor]{\hspace{1cm}{\sc Lemma}}
\newtheorem{conj}[theor]{\hspace{1cm}{\sc Conjecture}}
\newtheorem*{utver*}{\hspace{1cm}{\sc Proposition}}
\theoremstyle{definition}
\newtheorem{defin}[theor]{\hspace{1cm}{\sc Definition}}
\newtheorem{exa}[theor]{\hspace{1cm}{\sc Example}}
\newtheorem{rem}[theor]{\hspace{1cm}{\sc Remark}}
\newtheorem{prb}[theor]{\hspace{1cm}{\sc Problem}}

\newcommand{\Vol}{\mathop{\rm Vol}\nolimits}

\newcommand{\codim}{\mathop{\rm codim}\nolimits}

\newcommand{\rk}{\mathop{\rm rk}\nolimits}

\newcommand{\conv}{\mathop{\rm conv}\nolimits}

\newcommand{\Trop}{\mathop{\rm Trop}\nolimits}

\newcommand{\In}{\mathop{\rm in}\nolimits}

\def\R{\mathbb R}
\def\N{\mathbb N}
\def\Z{\mathbb Z}
\def\Q{\mathbb Q}
\def\C{\mathbb C}
\def\CC{({\mathbb C}\setminus 0)}
\def\TT{({\mathbb T}\setminus \mathfrak{0})}
\def\T{\mathbb T}

\def\CP{\mathbb C\mathbb P}
\def\RP{\mathbb R\mathbb P}

\begin{document}

\title{Characteristic classes of affine varieties and Pl\"ucker formulas for affine morphisms}
\author{Alexander Esterov\thanks{National Research University Higher School of Economics\newline This study (research grant No 14-01-0152) is supported by the National Research University--Higher School of Economics Academic Fund Program in 2014/2015. Partially supported by RFBR, grant 13-01-00755, MESRF, grant MK-6223.2012.1, and the Dynasty Foundation fellowship.}
}\date{}
\maketitle{}

\begin{abstract}
An enumerative problem on a variety $V$ is usually solved by reduction to intersection theory in the cohomology of a compactification of $V$. However, if the problem is invariant under a ``nice'' group action on $V$ (so that $V$ is spherical), then many authors suggested a better home for intersection theory: the direct limit of the cohomology rings of all equivariant compactifications of $V$. We call this limit the {\it affine cohomology} of $V$ and construct {\it affine characteristic classes} of 
subvarieties of a complex torus, taking values in the affine cohomology of 
the torus.  

This allows us to make the first steps in computing {\it affine Thom polynomials}. 
Classical Thom polynomials count how many fibers of a generic proper map of a smooth variety have a prescribed collection of singularities, and our affine version addresses the same question for generic polynomial maps of affine algebraic varieites.
This notion is also motivated by developing an intersection-theoretic approach to tropical correspondence theorems: they can be reduced to the computation of affine Thom polynomials, because 
the fundamental class of a variety in the affine cohomology is encoded by the tropical fan of this variety. 

The first concrete answer that we obtain is the affine version of what were, historically speaking, the first three Thom polynomials -- the Pl\"ucker formulas for the degree and the number of cusps and nodes of a projectively dual curve. This, in particular, characterizes toric varieties, whose projective dual is a hypersurface, computes the tropical fan of the variety of double tangent hyperplanes to a toric variety, 
and describes the Newton polytope of the hypersurface of non-Morse polynomials of a given degree. 
We also make a conjecture on the general form of affine Thom polynomials -- a key ingredient is the $n$-ary fan, generalizing the secondary polytope. \vspace{1ex}


\end{abstract}

\vspace{-5ex}
 
{\small
\tableofcontents}

\section{Introduction}

\subsection{Affine Pl\"ucker formulas}

Thom polynomials compute how many fibers of a generic proper morphism of a smooth variety have a prescribed collection of singularities. More generally, they count the cohomology class of a multisingularity stratum (the closure of the set of all points, whose preimages have a prescribed collection of singularities). Historically, the first Thom polynomials are the Pl\"ucker formulas, expressing the degree and the number of cusps and double points of the curve $\mathcal{C}^\vee$ projectively dual to a generic plane curve $\mathcal{C}$ in terms of the degree of $\mathcal{C}$. To see these expressions as Thom polynomials, note that $\mathcal{C}^\vee$ is the discriminant of  the tautological projection $$\pi:\{ (l,x)\in\CP^{2\vee}\times\CP^2\;|\; x\in l\cap\mathcal{C}\} \;\to\; \CP^{2\vee},\eqno{(*)}$$ i.e. the closure of all points $l\in\CP^{2\vee}$, whose fiber under $\pi$ has one singularity of local degree 2 ($A_1$ singularity). Similarly, the cusps and double points of $\mathcal{C}^\vee$ are the images of the fibers with one singularity of local degree 3 ($A_2$ singularity) or two $A_1$ singularities respectively. So the dual curve $\mathcal{C}^\vee$, the set of its cusps and the set of its double points are the multisingularity strata of the projection $\pi$.

The contemporary versions of the Pl\"ucker formulas are the Thom polynomials, expressing the fundamental classes of the  multisingularity strata for a generic map $\pi:M\to N$ of arbitrary compact smooth surfaces in terms of the characteristic classes of these surfaces: that is, the fundamental classes of the three multisingularity strata in $H(N)$ equal the direct images of $c_1,\, c_1^2+c_2$ and $-2(2c_1^2+c_2)$ respectively, where $c_i$ are the characteristic classes of the virtual vector bundle $\pi^*TN-TM$. We give \cite{K03} as a general reference for this fact and all subsequent appeals to Thom polynomials and multisingularity theory. In the special case of the tautological projection $(*)$, the characteristic classes are functions of the degree $d$ of the given curve $\mathcal{C}$, the fundamental class of the stratum $\mathcal{C}^\vee$ in $H_1(\CP^2;\Z)=\Z$ equals $d(d-1)$ and counts its degree, and the fundamental classes of the 0-dimensional 
multisingularity strata in $H_0(\CP^2;\Z)=\Z$ equal $3d(d-2)$ and $(d-3)(d-2)d(d+3)/2$ and count their cardinality.

We shall make the first steps towards computation of {\it affine Thom polynomials} that count how many fibers of a generic  polynomial map of an affine algebraic variety have a prescribed collection of singularities. An example of the affine version of the Pl\"ucker formulas is as follows.

\begin{exa} \label{exa1} 
Let $a:\{1,\ldots,n+1\}\to\Z$ be a non-negative concave function, $a(n+1)=0$, and $c_i$ be a generic univariate polynomial of degree $a(i)$. For how many values of the parameter $t\in\C$, is the polynomial $f_t(x)=c_1(t)x+\ldots+c_n(t)x^n+x^{n+1}$ not Morse? 

To answer this question, let $\pi$ be the restriction of the projection $\CC^3\to\CC^2,\, (x,y,t)\mapsto(x,t)$, to the graph $M=\{(x,y,t)\,|\, f_t(x)=y\}$, and let $A_1$, $A_2$ and $2A_1$ be its open multisingularity strata (i.e. the sets of points $(y,t)\in\CC^2$ such that the equation $f_t(x)=y$ has exactly one root of multiplicity 2, exactly one root of multiplicity 3, and exactly two roots of multiplicty 2 respectively). The sought number of non-Morse polynomials equals the cardinality $|2A_1|+|A_2|$, which can be found from the following three equations. Denote the Euler characteristics by $e$, then 

$$e(A_1)+2|2A_1|+2|A_2|=a(1)-2\sum_ia(i),$$

$$e(A_1)+2|2A_1|+|A_1|=3a(1)-4\sum_ia(i),$$

$$e(A_1)-|A_2|=a(1)-2\sum_i(3i-2)a(i).$$
\end{exa}

{\sc Proof.} Let $N\subset\R^3$ be the Newton polytope of the polynomial $F(x,y,t)=y-f_t(x)$, and let $D\subset\CC^2$ be the set of critical vaues of $\pi$ (i.e. the union of the strata $A_1$, $A_2$ and $2A_1$). Then, by the Kouchnirenko--Bernstein--Khovanskii formula \cite{khov0}, the Euler characteristic of the surface $F=0$ in $\CC^3$ equals the lattice volume of $N$. This is the first of the sought equations, because, counting $e(F=0)$ fiberwise, we obtain $e(F=0)=n\cdot e\bigl(\CC^2\setminus D\bigr)+(n-1)e(A_1)+(n-2)e(2A_1)+(n-2)e(A_2)=n\cdot e\CC^2-e(A_1)-2|2A_1|+2|A_2|=-e(A_1)-2|2A_1|+2|A_2|$ by the additivity of the Euler characteristic.

The set of critical points of $\pi$ is given by the equations $F=x\frac{\partial F}{\partial x}=0$. These equations are degenerate with respect to the Newton polytope $N$ in the sense of  \cite{khov0}, because the closure $C$ of this curve in the $N$-toric variety  intersects the 1-dimensional orbit, corresponding to the edge $[(1,0,0),\, (1,0,a(1))]$, at $a(1)$ points. However, one can easily check that, at every such point, $C$ transversaly intersects each of the adjacent 2-dimensional orbits. Thus, the Euler characteristic of  $F=x\frac{\partial F}{\partial x}=0$ is by $a(1)$ greater than the Euler characteristic of its perturbation $F=x\frac{\partial F}{\partial x}+G=0$, where $G$ is a small generic polynomial with Newton polytope $N$. The latter system of equations is non-degenerate with respect to $N$, so the Euler characteristic of its solutions equals $-2\Vol(N)$ by \cite{khov0}. This gives the second of the sought equations, because the Euler characteristic of the critical locus of $\pi$ equals $e(A_1)+2|2A_1|+|A_1|$ by the additivity of the Euler characteristic.

Finally, by \cite{EK}, the Newton polygon $N_D$ of $D$ is the fiber polytope of $N$. Explicitly, it has the vertices $(0,0)$, $(n,0)$ and $\bigl(i-1,2a(n)+2a(n-1)+\ldots+2a(i)\bigr)$ for $i=1,\ldots,n$. The closure of $D$ in the $N_D$-toric surface has three types of singularities: $|2A_1|$ transversal self-intersections, $|A_2|$ simple cusps and $a(1)$ simple tangencies with the 1-dimensional orbit $y=0$. Passing from $D$ to its generic perturbation $\tilde D$ with the same Newton polygon $N_D$, the singularities resolve and decrease the Euler characteristic of the curve by $|2A_1|$, $2|A_2|$ and $a(1)$ respectively. The resulting Euler characteristic of $\tilde D$ equals $-\Vol N_D$ by \cite{khov0}, which gives the third of the sought equations. $\quad\Box$
\vspace{1ex}

The statemnts regarding transversality and singularity types in this reasoning are justified by the results of Sections \ref{Scci} and \ref{Smain0}. The same computation is applicable to count the multisingularities of the projection $M\subset\CC^3\to\CC^2$ for a generic surface $M$ with an arbitrary Newton polytope $N$, although the answer is lengthy. 

\begin{rem} \label{rem00} 
Note that we could not reduce this example to classical Thom polynomials by compactifying the ambient tori $\CC^3\supset M$ and $\CC^2$: any such compactification would give a highly non-generic projection of the closure of $M$, to which classical Thom polynomials are not applicable. For example, if $M$ is a generic surface with Newton polytope $\conv\{(-2,0,0),\, (2,0,1),\, (0,1,0),\, (2k-1,1,1),\, (2k,1,1)\}$, and $D\subset\CC^2$ is the set of critical values of the projection $M\subset\CC^3_{x,y,t}\to\CC^2_{y,t}$, then the closure of $D$ in an arbitrary toric compactification of $\CC^2$ has a highly non-generic cusp with the Milnor number $2k-6$. 
The aim of the present paper is to develop a version of the characteristic classes and Thom polynomails applicable in such settings. 
 \end{rem}

\subsection{Affine cohomology ring}

In order to obtain an affine version of the Pl\"ucker formula for degree, we have to choose a cohomology-like ring, in which we shall compute the fundamental class of the discriminant of $\pi$. Remark \ref{rem00} suggests that the Chow ring of any particular compactification of $\CC^2$ is not a good choice from the perspective of classical Thom polynomials. For this reason we shall switch to the direct limit of the cohomology of all toric compactifications of a complex torus $\CC^n$. We shall call this graded ring the {\it affine cohomology} and denote it by $C$. The ring $C$ governs all equivariant enumerative problems in $\CC^n$ in the sense that the subvarieties $U$,  $V\subset\CC^n$ have the same fundamental class in $C$ if and only if every subvariety of complementary codimension, shifted by a generic element of the torus, has the same intersection number with $U$ as with $V$. 
This fact was first noticed in \cite{dCP2} for certain homogeneous spaces different from complex tori, but all of the arguments in this important paper are actually applicable to arbitrary spherical homogeneous spaces, including $\CC^n$ (see also \cite{brionduke89} and \cite{fultetal94} for details on the spherical setting).

The affine cohomology ring of the complex torus is isomorphic to the polytope algebra, which had been introduced in \cite{mcmullen89} and extensively studied (\cite{morelli93}, \cite{khpukhl93}) 
by the time when these two rings were identified in \cite{FS} and \cite{brion}. 
In particular, computation of the fundamental class of an algebraic hypersurface $H\subset\CC^n$ in the affine cohomology of $\CC^n$ 
amounts to computation of the Newton polytope of the defining equation of $H$.

\begin{exa} \label{exafb} Denoting the Newton polytope of $M$ from Example \ref{exa1} by $N$, the fundamental class of the discriminant $A_1$ (i.e. the Newton polygon of this plane curve) is the fiber polytope of $N$, see \cite{GP}. As explained above, one can see this equality as the affine version of the Pl\"ucker formula for degree in the setting of Example \ref{exa1}. 
\end{exa}

Recall that Minkowski summation of polytopes $A+B=\{a+b\,|\, a\in A, b\in B\}$ obviously extends to the notion of Riemann integral of a polyhedral-valued function, and the {\it fiber polytope} of $S\subset\R^3$ in a plane $L\subset\R^3$ is the Riemann integral $\int_{I} S(t)\, {\rm d}t$, where the segment $I\subset\R^3/L$ is the projection of $S$ along $L$, the volume form ${\rm d}t$ on the line $\R^3/L$ is induced by the lattice $\Z^3/L$, and $S(\cdot)$ is a polygon-valued function on $I$ that sends every $t$ to the plane section $S\cap(L+t)$ (see \cite{bs} for details).

\begin{rem} More generally, computation of the affine fundamental class of an arbitrary subvariety $H\subset\CC^n$ in the affine cohomology of $\CC^n$ amounts to computation of the {\it tropical fan} of $H$ (see \cite{Kaz99} and \cite{Kaz}, \cite{ST}, \cite{markw07}). However, we prefer to think here in terms of affine cohomology classes, because some of our constructions and reasonings work well in the context of arbitrary spherical spaces, to which the notion of tropical fan has not been generalized so far.
\end{rem}

\subsection{Affine multisingularity theory}

The setting in Example \ref{exa1} admits many variations: we can consider more complicated $M\subset\CC^3$, increase the dimensions of the complex tori, replace the projection of tori with a generic polynomial map, look for the fundamental classes of more complicated singularity strata and so on. However, all of these variations can be reduced to the universal Problem \ref{prb1} which follows. The reduction is outlined in Section \ref{Sred} and will be treated in detail in a separate paper. 

Choose a finite set $A$ in the character lattice $\Z^n$ of the torus $\CC^n$, and consider the space $\C^A$ of all linear combinations of the characters from $A$. 
For every finite collection $S$ of isolated hypersurface singularities, consider the {\it universal $S$-multisingularity stratum} $\{S\}\subset\C^A$, i.e. the set of all $f\in\C^A$, such that the singularities of the hypersurface $f=0$ in $\CC^n$ are in one-to-one correspondence with the singularities of $S$ and are equivalent to them; see Section \ref{Stmt} for a more precise definition. 
\begin{prb}\label{prb1}
Compute the fundamental class of the universal $S$-multisingularity stratum $\{S\}\subset\C^A$ in the affine cohomology of the torus $\CC^A$. 
\end{prb}
Many recent works can be seen as solutions of this problem for different $A$ and $S$, see Section \ref{Stmt} for an overview. We now describe our contribution and its immediate applications.

Corollary \ref{newtdiscr} solves Problem \ref{prb1} for the universal codimension 1 stratum: it describes the Newton polytope of the $A$-discriminant. Although this Newton polytope has been extensively studied (\cite{GKZ94}, \cite{GP}, \cite{DFS}, \cite{Edcg}, \cite{Smix}, \cite{Dmix}, \cite{Tsikh}, etc.), even in this special case our approach yields something new: Corollary \ref{newtdiscr} is a positive formula for the tropicalization of the $A$-discriminant, different from the known one \cite{DFS} 
and sufficient to classify non dual defective toric varieties.

Recall that a projective variety is called dual defective, if its projectively dual variety is not a hypersurface. Classification of dual defective varieties is a classical problem, and even the toric case remains unsolved (smooth dual defective toric varieties are classified in \cite{sdr}, see also \cite{cc} and \cite{sdr2} for further results). We classify non dual defective toric varieties. A set $B\subset\Z^n$ is called a {\it circuit}, if it consists of $n+2$ points, any $n+1$ of whom are affinely independent. It is called an {\it iterated circuit}, if there exists a flag $$\varnothing=L_0\subset\{0\}=L_1\subset L_2\subset\ldots\subset L_k=\Z^n$$ such that, for every $i=1,\ldots,k$, the set $B\cap L_i\setminus L_{i-1}$ consists of $\dim L_i-\dim L_{i-1}+1$ points, whose images under the projection $L_i\to L_i/L_{i-1}$ together with $0$ form a circuit.
\begin{theor} (Corollary \ref{newtdiscr}.2) An $A$-toric variety is not dual defective if and only if, up to a translation, $A$ contains an iterated circuit.
\end{theor}

Problem \ref{prb1} for the universal codimension 2 strata $A_2$ and $2A_1$ is solved by Theorems \ref{th1}, \ref{th2} and \ref{th3}, which generalize the three equations from Example \ref{exa1}. The answer is illustrated for a particular simple $A\subset\Z^2$ in Section \ref{Sexample}, see also \cite{katz09} for another example. 

In particular, this allows the description of the Newton polytope of the hypersurface $H\subset\CC^A$ of all polynomials $f:\CC^n\to\C$ that are not Morse functions (i.e. have two equal critical values or a degenerate critical value), see  Section \ref{Sred}. At the level of degrees (i.e. when  $A=d\cdot($standard simplex$)$, and we are computing the degrees of the sets $A_2$, $2A_1$ and $H$ in terms of $d$), these problems were solved in \cite{aluf98} and \cite{landozvon} (for $n=1$) respectively. 

A recent paper \cite{DHT} suggests a different approach to the computation of the fundamental class of the stratum $2A_1$ for $n=1$ (i.e. the set of all univariate polynomials with two multiple roots), and greatly clarifies the answer in this case.

\subsection{Affine characteristic classes}

Multisingularity theory is based on characteristic classes, and the proof of the aformentioned results will be based on
{\it affine characteristic classes} of subvarieties of a complex torus, taking values in the affine cohomology ring $C$. 
\begin{theor} \label{th00} For every algebraic set $V\subset\CC^n$, the Schwarz-MacPherson characteristic classes of $V$ in the toric compactifications of $\CC^n$ converge to a certain class $\alpha\in C$: there exists a toric compactification $X_0\supset\CC^n$ such that the restriction of $\alpha$ to the cohomology of every richer smooth toric compactification $X\mapsto X_0$ is Poincare dual to the Schwarz-MacPherson characteristic class of $V$ in $X$. 
\end{theor}
See the original papers \cite{mcph} and \cite{sch82} and e.g. a survey \cite{alufC} for the background on Schwarz-MacPherson classes.
\begin{defin} The limit in Theorem \ref{th00} is called the {\it affine characteristic class} of $V$ and will be denoted by $\langle V\rangle$.
\end{defin}
This class inherits nice functorial properties of the Schwarz-MacPherson classes: e.g. $\langle U\cap V\rangle+\langle U\cup V\rangle=\langle U\rangle+\langle V\rangle$ and $\langle U\cap gV\rangle=\langle U\rangle\langle V\rangle$ for a generic $g\in\CC^n$. Denoting the codimension $i$ component of $\langle V\rangle$ by $\langle V\rangle_i$, we have $\langle V\rangle_n=e(V),\; \langle V\rangle_i=0$ for $i<\codim V$, and $\langle V\rangle_{\codim V}$ is the fundamental class of $V$ in $C$.
\begin{exa} The formal rational function $\frac{N_1\ldots N_k}{(N_1+1)\ldots(N_k+1)}$ of the Newton polytopes of the polynomials $f_1,\ldots,f_k$, which was considered in \cite{khov0} and in subsequent literature in relation to the nondegenerate complete intersection $f_1=\ldots=f_k=0$ in a complex torus, turns out to equal the affine characteristic class of this complete intersection (when evaluated in the affine cohomology ring), see Examples \ref{trophyp} and \ref{kbkex} for more details.
\end{exa}

\begin{rem} 
Affine characteristic classes seem to be indispensable in solving Problem \ref{prb1}, because they appear of their own accord even in the answer for the universal strata $A_2$ and $2A_1$ of codimension 2. More specifically, the answer given by Theorems \ref{th1}, \ref{th2} and \ref{th3} is a system of three independent linear equations on the fundamental classes of $A_2$ and $2A_1$ and one more unknown class, which turns out to be the second affine characteristic class of the discriminant $A_1$. Note the similarity to Example \ref{exa1}, where the third unknown $e(A_1)$ can also be interpreted as the second affine characteristic class of the discriminant.
\end{rem}
\begin{rem} Note that, in Theorem \ref{th00}, we take the Schwarz-MacPherson classes of $V$ itself, and not the closure of $V$ in $X$. The statement of the theorem obviously generalizes to an arbitrary spherical homogeneous space, however I have no proof. For the toric case, the proof occupies Section \ref{Sconstr1} and is indirect: firstly, affine characteristic classes are constructed in a different way and shown to have functorial properties, similar to Schwarz-MacPherson classes. Secondly, these properties are shown to imply that the restriction of the affine characteristic class to the cohomology ring of a sufficiently rich compactification equals the Schwarz-MacPherson class.
\end{rem}

\subsection{Relation to tropical enumeration}

The classical approach to enumerative geometry is through intersection theory, characteristic classes and Thom polynomials. Answers to enumerative questions can usually be given by certain Thom-like polynomials, see e.g. \cite{G98} and \cite{T10} for the problem of counting rational curves passing through given generic points of a given surface. 
In the last decade, a fundamentally new class of enumerative results was invented by Mikhalkin in \cite{M}. Results of this type are seemingly unrelated to the aforementioned classical approach to enumerative geometry, and are referred to as tropical correspondence theorems. Such theorems state that the answers to certain enumerative questions coincide over $\C$ and over the tropical semifield $\T=(\R\cup\{-\infty\}, \max, +)$. Since enumerative geometry over $\T$ is combinatorics, this gives a new combinatorial answer to classical problems of enumerative geometry over $\C$.

Proofs of tropical correspondence theorems usually rely upon some kind of deformation theory technique. Every such proof, besides establishing the desired equality between the number of complex and tropical objects, implicitly contains an algorithm of how to reconstruct the whole 1-parametric family of complex objects being deformed to a given tropical one. E. g. for a given 1-parametric family of points $P(t)$ in the projective plane and a given rational tropical curve $T$, passing through the points of the tropicalization of $P(t)$, the deformation-theoretic proof of Mikhalkin's correspondence theorem  (\cite{Sh}, \cite{sn}, \cite{tyomkin}) provides an algorithm of how to compute (up to an arbitrary power of $t$) the 1-parametric family of rational curves, passing through $P(t)$ and tending to $T$.
One might hope that another approach to the proof, which does not inevitably provide such excessive information, would appear to be more universally applicable.
 
One such intersection theoretic approach was recently developed in \cite{grossR} and references therein. It is oriented towards enumeration of rational curves.

We suggest yet another approach, oriented towards curves and hypersurfaces, given by implicit equations. It turns out that once we know the fundamental class of a universal $S$-multisingularity stratum in the affine cohomology, we can obtain tropical correspondence theorems for hypersurfaces with multisingularity $S$ and various incidence conditions. So tropical enumeration can be regarded as the classical approach to enumerative geometry, in which the cohomology, the characteristic classes and the Thom polynomials are replaced by their affine versions. This is because the affine fundamental class of a subvariety of a torus is encoded by its {\it tropicalization} (\cite{Kaz99} and \cite{Kaz}, \cite{ST}, \cite{markw07}). 
See Section \ref{Senum} for details, and a simple illustration below:
we deduce Mikhalkin's correspondence theorem for curves with one node from the description of the Newton polyhedron of the A-discriminant (i.e. the affine fundamental class of the universal singularity stratum of codimension 1). We start with the definition of tropical numbers. 

Let $\T=\R\sqcup\{-\infty\}$ be the tropical semifield with the operations \begin{center} $a\underset{\T}{\cdot}b=a\underset{\R}{+} b\quad$ and $\displaystyle\quad a\underset{\T}{+}b=\left\{ {\max(a,b) \;\; \mbox{ if } a\ne b \atop [-\infty,\, a] \;\;\;\; \mbox{ if } a=b}\right.$ \end{center}

\begin{rem} Multivalued tropical summation that we use is a little more convenient than the conventional one (see e.g. the subsequent definition) and is also more natural: the tropical summation and multiplication should simulate the behaviour of the degree of the sum and the product of two polynomials, whereas $\deg(f+g)$ can attain any value within $[-\infty,\deg f]$ if $\deg f = \deg g$.\end{rem}

In what follows, $\mathfrak{1}\in\T$ always stands for $0\in\R$, and $\mathfrak{0}\in\T$ for $-\infty$. Moreover, due to multivaluedness, ``$=\mathfrak{0}$'' in the tropical context always means ``$\ni-\infty$''.

\begin{rem} We shall denote tropical objects in Fraktur, because we think of the tropical torus $\T\setminus\mathfrak{0}$ as the Lie algebra of the complex torus $\C\setminus 0$.
\end{rem}
\begin{defin} A {\it tropical hypersurface} is a set of the form $S = \{ x \,|\, f(x){=}\mathfrak{0}\}$, where $f$ is a tropical polynomial. In a small neighborhood of a generic point $x\in S$, the polynomial $f$ can be represented as $g\cdot h^k$, where $g(x)\ne \mathfrak{0}$ and $h$ is irreducible. The integer $k$ is called the {\it multiplicity} of the hypersurface $S$ at the point $x$.
\end{defin}

For example, the following picture shows the graph of $f(x)=x^3+x^2+6x+8$ in bold lines, and the set $\{f=\mathfrak{0}\}=\{2,\, 3\}$ with the multiplicities $1$ and $2$ respectively as white dots.
\begin{center}
\noindent\includegraphics[width=5cm]{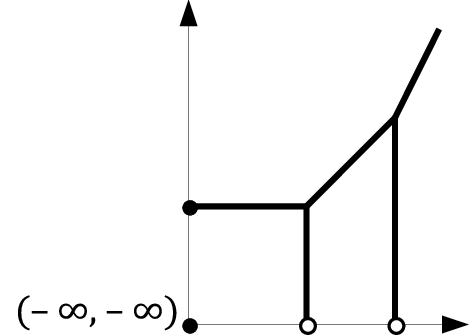}
\end{center}

Tropical hypersurfaces $H_1,\ldots,H_n$ in the tropical torus $\R^n=(\T\setminus 0)^n$ are said to intersect {\it transversally}, if they intersect at finitely many points, and, in a small neighborhood of every such point $p$, every $H_i$ is given by a tropical equation of the form $x^{a_{p,i}}+p^{a_{p,i}}=\mathfrak{0}$ for some $a_{p,i}\in\Z^n$. The {\it intersection number} of $H_1,\ldots,H_n$ is then defined as the sum $\sum_p |\det(a_{p,1},\ldots,a_{p,n})|$. The tropical Kouchnirenko--Bernstein formula states that the intersection number equals the mixed volume of the Newton polytopes of $H_1,\ldots,H_n$.

\begin{exa}[Mikhalkin's correspondence for curves with one node]
For a given $A\subset\Z^2,\, 0\in A$, we wish to count polynomials $f\in \C^A$ such that the curve $f=0$ has one singularity and passes through a given collection of generic $q=|A|-2$ points $p_1,\ldots,p_q\in\CC^2$. In other words, we wish to count the intersection number $I$ of the following hypersurfaces in $\C^A$:

the incidence conditions $H_1,\ldots,H_q$, where $H_i=\{f \,|\, f(p_i)=0\}$;

the normalization $H_0=\{f \,|\, f(0)=1\}$;

the $A$-discriminant $S=\{f \,|\, f=0$ is not regular$\}$.

Let us tropicalize these objects: choose points $\mathfrak{p}_1,\ldots,\mathfrak{p}_q\in (\T\setminus 0)^2$ and define the tropical hypersurfaces $\mathfrak{H}_i=\{f \,|\, f(\mathfrak{p}_i)=\mathfrak{0}\}$ in $\T^A$ and $\mathfrak{H}_0=\{f \,|\, f(\mathfrak{0})=\mathfrak{1}\}$. Choose the tropical polynomial $\mathfrak{D}$ with unit coefficients and the same Newton polytope as the $A$-discriminant $S$, and define the tropical discriminant $\mathfrak{S}$ by the 
equation $\mathfrak{D}=\mathfrak{0}$. For generic choice of ${\mathfrak p}_1,\dots, {\mathfrak p}_q$, the tropical hypersurfaces $\mathfrak{H}_0,\ldots,\mathfrak{H}_q,\mathfrak{S}$ are transversal, and we denote their intersection number by $\mathfrak{I}$. Then we have
$$I=\mathfrak{I},$$
because both parts equal the mixed volume of the Newton polytopes of $H_0,\ldots,H_q$ and $S$ by the Kouchnirenko--Bernstein formula over $\C$ and $\T$ respectively. 

The description of $\mathfrak{S}$ (see Section \ref{Sdiscr}) shows that the intersection points $\mathfrak{f}\in \mathfrak{H}_0\cap\ldots\cap \mathfrak{H}_q\cap \mathfrak{S}$ are exactly the equations of the tropical curves with one node in the sense of Mikhalkin, passing through $\mathfrak{p}_1,\ldots,\mathfrak{p}_q$. Moreover, the intersection multiplicity of $\mathfrak{H}_0,\ldots,\mathfrak{H}_q$ and $\mathfrak{S}$ at every such intersection point $\mathfrak{f}$ equals the multiplicity, assigned to the curve $\mathfrak{f}=\mathfrak{0}$ by Mikhalkin's correspondence theorem.
\end{exa}

{\sc Acknowledgements.} Many steps in this work are inspired by collaboration with A. Khovanskii: in particular, I learned from him the aformentioned interpretation of tropical correspondence theorems. 
The codimension 1 part of the equality in Proposition \ref{eqstar} is essentially the key observation in our join work \cite{EK}. The relation of tropical characteristic classes to Schwartz-MacPherson classes, conjectured in the first version of this paper, turned into a new simple proof of key Proposition \ref{schoncc} thanks to B. Sturmfels, who suggested the connection of this work with \cite{huh}. The natural generality for spherical characteristic classes was suggested by V. Kiritchenko. Attention to this work from F. Block, A. Dickenstein, S. di Rocco, B. Nill and L. Tabera has greatly improved the first verison of the paper.

\section{Affine characteristic classes} \label{SSacc}

Affine intersection theory and affine characteristic classes are introduced for arbitrary spherical homogeneous spaces in Section \ref{Sscc} and specialized to the toric case in Sections \ref{Stit} and \ref{Stcc} respectively. The rest of Section \ref{SSacc} is more technical: Sections \ref{Spfp}---\ref{Sconstr1} are devoted to the proof of existence of affine characteristic classes in the toric setting, and Section \ref{Scci} is devoted to characteristic classes of degenerate complete intersections of a certain type, whose importance will become clear in Section 3.

\subsection{Affine intersection theory and characteristic classes} \label{Sscc}

Let $G$ be a reductive algebraic group (i. e. the complexification of a real compact connected group), and let $X$ be an $n$-dimensional spherical homogeneous $G$-space (i. e. a space with a transitive action of $G$, such that some Borel subgroup of $G$ has a dense orbit in $X$). For algebraic subsets $P$ and $Q\subset X$ of complementary dimensions $p+q=n$, the set $gQ=\{gq\,|\, q\in Q\}$, the shifting of $Q$ by a generic element $g\in G$, meets $P$ at finitely many points, and we denote the number of these points by $|P\cdot Q|$.

Denote the space of linear combinations of irreducible codimension $k$ algebraic subsets of $X$ by $Z^k$; then the aforementioned pairing extends to a pairing $|\cdot|:Z^p\times Z^q\to\Z$ by linearity, and the space $C=\bigoplus_k C_k=\bigoplus_k Z_k/\{P\, |\, \forall Q\; |P\cdot Q|=0\}$ has a natural ring structure, compatible with this pairing (see \cite{dCP2},\cite{brionduke89}, \cite{fultetal94}): for algebraic subsets $R$ and $S\subset X$ of arbitrary dimensions, the intersection of $R$ and $gS$ represents the same element of $C$ for almost all $g\in G$, and this element is called the {\it product} $R\cdot S$. 
We shall call this ring the {\bf affine cohomology}, or, following \cite{dCP2}, the {\bf ring of conditions} of $X$.

\begin{rem} Recall that this multiplication does not exist for non-spherical spaces: for example, if  $G=X=\C^3$, then two lines represent the same class in $C$ if and only if they are parallel, however the intersection of the surface $z=xy$ and the plane $x=0$, shifted by $g\in\R^3$, is a line, whose direction depends on $g$.
\end{rem}

\begin{defin} \label{defspher} The {\bf affine characteristic class} is a mapping that sends every algebraic subset $V\subset X$ to an element $\langle V\rangle=\langle V\rangle_0+\ldots+\langle V\rangle_n\in C,\, \langle V\rangle_i\in C_i$, with the following properties:
\end{defin}

(1) If $V\subset X$ has codimension $k$, then $\langle V\rangle_i=0$ for $i<k$, $\langle V\rangle_k$ is the equivalence class of $V$ in $C_k$, and $\langle V\rangle_n\in C_0=\Z$ is the Euler characteristic $e(V)$.

\vspace{1ex}

(2) For any $U$ and $V\subset X$ and generic $g\in G$, we have $\langle U\cap gV\rangle=\langle U\rangle\langle V\rangle/\langle X\rangle$. Note that $\langle X\rangle$ is invertible by property (1).

\vspace{1ex}

(3) For any spherical homogeneous spaces $X$ and $Y$ and any algebraic subsets $U\subset X$ and $V\subset Y$, we have $\langle U\times V\rangle=\langle U\rangle\times \langle V\rangle$.

\vspace{1ex}

(4) The mapping that assigns $\langle V\rangle$ to the characteristic function of $V$, extends by linearity to the space of all constructible functions $X\to\Z$, i.e. $\langle U\cap V\rangle+\langle U\cup V\rangle=\langle U\rangle+\langle V\rangle$. Recall that a {\it constructible function} is a linear combination of the characteristic functions of algebraic sets.

\vspace{1ex}

(5) For a morphism $p:X\to Y$ of spherical homogeneous spaces and an algebraic subset $V\subset X$, we have $p_*\langle V\rangle=\langle p_* V\rangle$. Here $p_* V:Y\to\Z$ is the {\it MacPherson direct image} of $V$, whose value at every $y\in Y$ is defined as $e\Bigl(p^{-1}(y)\cap V\Bigr)$.

\vspace{1ex}
(6) For a smooth equivariant compactification $\overline{X}\supset X$, such that the affine characteristic class $\langle V\rangle$ is contained in the cohomology $H^\bullet(\overline{X})\subset C$, this class is Poincare dual to the {\it Schwarz-MacPherson class} of $V$ in $\overline{X}$, see \cite{mcph} and \cite{sch82}. Recall that the natural inclusions $H^\bullet(\overline{X})\subset C$ are induced by the fact that $C$ is the direct limit of $H^\bullet(\overline{X})$ over all equivariant compactifications $\overline{X}\supset X$.

\vspace{1ex}

Note that the affine characteristic class is uniquely defined by property (6) and by properties (1--5).
\begin{theor} \label{thm23} In the toric case $G=X=\CC^n$, the affine characteristic class exists, and property (2) reads as $\langle U\cap gV\rangle=\langle U\rangle\langle V\rangle$, because $\langle X\rangle=1$.
\end{theor}
The proof is given in Section \ref{Sconstr1} and is based on an explicit combinatorial model for the ring of conditions $C$ of a complex torus. This combinatorial model is widely known as the {\bf ring of tropical fans} (\cite{ST}, \cite{markw07}), although they are essentially the same thing as Minkowski weights from \cite{FS} and c-fans from \cite{Kaz}. The affine characteristic class, considered as an element of the ring of tropical fans, will be referred to as the {\bf tropical characteristic class}.

\begin{rem} 
I cannot generalize Theorem \ref{thm23} to an arbitrary spherical space $X$, but this formal problem is not the main difficulty with affine characteristic classes for arbitrary $X$. For instance, the class $\langle X\rangle$, which is especially important due to property (2), has already been cosidered in \cite{valya} under the name of the {\it non-compact characteristic class of} $X$ (see also \cite{brion2} and \cite{brion3} for its equivariant version), but it is not yet computed even for $X=SL_n$ with large $n$.

\end{rem}

So, from now on, we restrict our attention to the toric case, recalling the notion of a tropical fan (we actually need the slightly more general notion of tropical fan with polynomial weights), and develop the technique for computation of tropical characteristic classes to the extent that we need in this paper.

\subsection{Toric intersection theory and tropical fans} \label{Stit}

The complex torus $G=X=\CC^n$ is the only spherical homogeneous space for which the structure of the ring of conditions is completely understood. More specifically, the definition of the intersection product as the class of $U\cap gV$ for generic $g\in G$ is not constructive because of the word `generic'. However, for $G=X=\CC^n$, we know how to check  explicitly if 
the intersection number of algebraic subsets $U$ and $V\subset\CC^n$ of complementary dimension 
equals the product of the classes of $U$ and $V$ in the ring of conditions. In what follows, we write this equality as 
$\langle U\cap V\rangle=\langle U\rangle\langle V\rangle$ in accordance with our notation for characteristic classes.

Let $L$ be the character lattice of the complex torus, and let $I$, $J\subset\C[L]$ be the radical ideals, defining $U$ and $V$. For a linear function $\gamma:L\to\Z$, define $\In_{\gamma} I$ as the ideal, generated by the initial terms of the polynomials from $I$ in the sense of the partial ordering $\gamma$ on the lattice of monomials. This ideal $\In_{\gamma} I$ defines a 
variety denoted by $\In_{\gamma} U$.
\begin{utver}[\cite{Kaz}] \label{tropint0} If $\In_{\gamma} U\cap \In_{\gamma} V$ is empty 
for every non-zero $\gamma\in L^*$, then $\langle U\cap V\rangle=\langle U\rangle\langle V\rangle$.
\end{utver}
In particular, this leads to a combinatorial representation of $C$ as the space of certain fans, based on the fact that the set $\{\gamma\, |\, \In_{\gamma} U\ne\varnothing\}$ is a polyhedral fan. We recall this combinatorial representation in the form of \cite{Emmj12}.

Let $L$ be an $n$-dimensional integer lattice. Consider a pair $(P,\varphi)$, where 

1) $P\subset L\otimes\R$ is a union of finitely many disjoint rational convex relatively open codimension $k$ polyhedral cones, so that $P$ coincides with a codimension $k$ plane in a neighborhood of any of its points $p$ (this plane is denoted by $T_p P$); 

2) The function $\varphi:P\to\R$ equals a rational polynomial $\varphi_p:T_p P\to \R$ in a small neighborhood of every $p\in P$.
 
Two such pairs $(P,\varphi)$ and $(Q,\psi)$ are said to be {\bf equivalent}, if $\varphi(p)=\psi(p)$ for every $p\in P\cap Q$ such that $T_p P=T_p Q$, and $\varphi=0$ on $P\setminus\overline{Q}$ and $\psi=0$ on $Q\setminus\overline{P}$.

For a rational $k$-dimensional plane $R\subset L\otimes\R$ and a point $x\in L\otimes\R$, the affine plane $R+x$  is said to be {\bf transversal} to $(P,\varphi)$, if it meets $P$ at finitely many points and does not meet $\overline{P}\setminus P$. The {\bf tropical intersection number} $(R+x)\cdot(P,\varphi)$ is then defined as
$$\sum_{p\in (R+x)\cap P} \varphi(p) \left| \frac{L}{(R\cap L)+(T_pP\cap L)} \right|.$$ 
This intersection number is a locally polynomial function of $x$, defined on an open dense subset of $L\otimes\R$. If this function extends to a continuous function $i_R:L\otimes\R\to\R$ for every $R$, then the pair $(P,\varphi)$ is said to be {\bf tropical}, and the {\bf tropical intersection number} $(R+x)\cdot(P,\varphi)$ is defined to be $i_R(x)$ for every $x$ by continuity, even if $R+x$ is not transversal to $(P,\varphi)$.

Equivalence classes of codimension $k$ tropical pairs are called tropical fans and form the set $\mathcal{K}_k(L)$. For any two fans $\mathcal P$ and $\mathcal Q$, there exists a unique (up to equivalence) fan $\mathcal S$, such that ${\mathcal  P}\cdot R+{\mathcal  Q}\cdot R={\mathcal  S}\cdot R$ for every affine plane $R$. This fan $\mathcal S$ is called the sum ${\mathcal P}+{\mathcal  Q}$, see \cite{Emmj12} for a more constructive definition. 

The set $\mathcal{K}_k(L)$ is a $\Q$-vector space with respect to this summation. It splits into the sum $\bigoplus_d \mathcal{K}_k^d(L)$, where $\mathcal{K}_k^d(L)$ consists of pairs $(P,\varphi)\in \mathcal{K}_k(L)$ such that $\varphi$ is locally a homogeneous polynomial of degree $d$. The direct sum $\bigoplus_k \mathcal{K}_k(L)$ is denoted by $\mathcal{K}(L)$ and will be called the {\bf space of tropical fans with polynomial weights}.

For tropical pairs $(P,\varphi)\in \mathcal{K}(L)$ and $(Q,\psi)\in \mathcal{K}(M)$, their {\bf Cartesian product} is defined as $(P\times Q,\varphi+\psi)\in \mathcal{K}(L\oplus M)$.
For an epimorphism of lattices $f:L\to M$ of dimension $n$ and $m$ respectively, the {\bf push-forward} $f_*({\mathcal P})$ of a tropical fan ${\mathcal P}\in \mathcal{K}_k(L),\, k\geqslant n-m$, is defined as the unique ${\mathcal S}\in \mathcal{K}_{k-n+m}(M)$ such that ${\mathcal P}\cdot f^{-1}(R)={\mathcal S}\cdot R$ for every affine $(k-n+m)$-dimensional plane $R\subset M\otimes\Q$. A more explicit description is as follows (see \cite{Emmj12}).
\begin{utver} \label{projtrop2} 
The image of the tropical fan $\mathcal{P}$ that consists of $(n-k)$-dimensional cones $C_i\subset L\otimes\R$ with multiplicities 
$m_i:C_i\to\R$ is a tropical fan $f_*\mathcal{P}$ that consists of all $(n-k)$-dimensional cones $f(C_i)$ with multiplicities $(m_i\circ f^{-1})\cdot\bigl|L\bigl/\bigr.(C_i+\ker f)\bigr|$. 
\end{utver}

The vector space $\mathcal{K}(L)$ has the natural differential ring structure: the {\bf intersection product} of ${\mathcal P}\in \mathcal{K}_k(L)$ and ${\mathcal Q}\in \mathcal{K}_m(L)$ is the unique ${\mathcal S}\in \mathcal{K}_{k+m}(L)$, such that ${\mathcal S}\cdot R=({\mathcal P\times \mathcal Q}\times R)\cdot($diagonal) in $(L\oplus L\oplus L)\otimes\Q$ for every $(k+m)$-dimensional affine plane $R\subset L\otimes\Q$. The {\bf corner locus derivation} $\delta:\mathcal{K}_k^d(L)\to \mathcal{K}^{d-1}_{k+1}(L)$ is the unique derivation such that $\delta({\mathcal P})\in \mathcal{K}^0_1(L)$ is the corner locus of a continuous piecewise linear function ${\mathcal P}\in \mathcal{K}^1_0(L)$, see \cite{Emmj12} for more constructive definitions.

Denote the subrings $\bigoplus_k \mathcal{K}_k^0(L)$ and $\bigoplus_d \mathcal{K}^d_0(L)$ of $\mathcal{K}(L)$ by $\mathcal{K}^0(L)$ and $\mathcal{K}_0(L)$ respectively. The former one is also known as the {\bf ring of tropical fans} (see e.g. \cite{FS}, \cite{Kaz}, \cite{ST}, \cite{markw07} for a more explicit description of the sum and the product in this ring), and the latter one is the ring of all continuous piecewise polynomial functions on $L\otimes\Q$ (with the conventional sum and product).
The subrings $\mathcal{K}_0(L^*)$ and $\mathcal{K}^0(L^*)$ provide well known combinatorial models for the ring of conditions $C$ of the complex torus $\CC^n$ with the character lattice $L$: $$C=\mathcal{K}^0(L^*)=\mathcal{K}_0(L^*)/I,$$ where $I$ is the ideal, generated by linear functions on $\Q^n$. The first of these models was found in \cite{FS}, \cite{mcm} and \cite{Kaz}, and the second one in \cite{stanley} and \cite{brion}. The isomorphism between the two models is established by the maps $\delta^k:\mathcal{K}^k_0(L^*)\to \mathcal{K}^0_k(L^*)$
(see also \cite{brion97} and \cite{paynekatz} for other descriptions of this isomorphism).
The element of $\mathcal{K}^0(L^*)$, corresponding to the class of a subvariety $V\subset\CC^n$ in the ring of conditions $C$, is called the {\bf tropical fan} $\Trop V$, and admits the following description:

\begin{utver}[\cite{Kaz}] \label{tropn}
Let $V$ be given by a radical ideal $I\in \C[L]$, then its tropical fan is represented by the pair $(P,\varphi)\in \mathcal{K}^0\bigl(L^*\bigr),\,P\subset\Q^n,\,\varphi:P\to\Q$, where $P=\{\gamma\, |\, \dim_\C\C[L]/\In_{\gamma} I$ is finite and positive$\}$, and $\varphi(\gamma)=\dim_\C\C[L]/\In_{\gamma} I$.
\end{utver}

\begin{exa} \label{trophyp}
If $N$ is the Newton polytope of a Laurent polynomial on the complex torus $\CC^n$, then the tropical fan of the hypersurface $f=0$ equals the corner locus of the support function of $N$. 
We denote this fan by $[N]\in\mathcal{K}^0_1(\Z^n)$ and call it the {\bf dual fan} of $N$.
\end{exa}

\begin{rem} \label{rcondfunct} This isomorphism between $\mathcal{K}^0(L^*)$ and $C$ respects epimorphisms of tori (\cite{Kaz}, \cite{ST}): it sends the epimorphism of the rings of conditions, induced by an epimorphism $\CC^n\to\CC^m$, to the epimorphism $\mathcal{K}^0(L^*)\to \mathcal{K}^0(M^*)$, induced by the corresponding epimorphism of the dual character lattices $L\to M$. \end{rem}

\subsection{Tropical characteristic classes} \label{Stcc}

\begin{defin} \label{tropcc}
The element of the ring of tropical fans $\mathcal{K}^0\bigl(L^*\bigr)$, corresponding to the affine characteristic class of a subvariety $V\subset\CC^n$, is called the {\bf tropical characteristic class} of $V$.
\end{defin}

\begin{rem} Here are some natural interesting questions about tropical characteristic classes that we do not address in this paper:

1. Construct ``higher Newton polytopes'' of a codimension $k$ variety $V$, i.e. codimension 1 tropical fans $B_1,B_2,\ldots$ such that $\langle V\rangle_{k+i}=\langle V\rangle_kB_1\ldots B_i$. 
Are they dual to convex polytopes if $V$ is smooth?

2. Express $\langle V\rangle$ in terms of an arbitrary (not necessarily toric!) compactification $C\supset\CC^n$, such that the closure of $V$ is transversal to the strata of $C$.

3. Define and study the ``generalized Severi variety'', whose points parameterize all subvarieties $V\subset\CC^n$ with a given tropical characteristic class. Classical Severi varieties are such spaces for plane curves. See \cite{kp11} for the first step in this direction at the level of fundamental classes.

4. Obtain a linear-algebraic description for the tropical characteristic class of a variety in terms of its defining ideal, similarly to Proposition \ref{tropn}.

5. Try to apply tropical characteristic classes to the conjectures in \cite{huh}.
\end{rem}

We now compute $\langle V\rangle$ under the assumption that $V$ is {\bf sch\"{o}n}, i.e. there exists a fan $\Sigma$ such that the closure of $V$ in the corresponding toric compactification $X_\Sigma\supset\CC^n$ is smooth and intersects the orbits of $X_\Sigma$ transversally (see for instance \cite{khovtor} and \cite{danil} for background on toric varieties and their relation to Newton polyhedra).
Note that our definition of Sch\"{o}n varieties is a bit more restrictive than in the original paper \cite{Tev}.

 For every cone $\Gamma\in\Sigma$, denote the intersection of the closure of $V$ with $\Gamma$-orbit of $X_\Sigma$ by $V_\Gamma$.
\begin{utver} \label{schoncc} 
If $V$ is sch\"{o}n, then the class $\langle V\rangle_i\in \mathcal{K}^0_i(L^*)$ is represented by the pair $(P,\varphi),\, P\subset\Q^n,\, \varphi:P\to\Q$, such that $P$ is the union of all codimension $i$ cones in $\Sigma$, and the value of $\varphi$ at every such cone $\Gamma$ is the Euler characteristic $e(V_\Gamma)$.
\end{utver}
The proof will be given in Section \ref{Sconstr1}, because it comes as a byproduct when we prove the existence of affine characteristic classes.

\begin{exa}\label{kbkex} In particular, if $V$ is a generic hypersurface with Newton polytope $\Delta$, then, counting $e(V_\Gamma)$ by the Kouchnirenko formula, we obtain $\langle V\rangle=\frac{[\Delta]}{1+[\Delta]}$, or $\langle V\rangle_i=-(-[\Delta])^i$. Recall that $[B]$ is the {\bf dual fan} of the polytope $B$, see Example \ref{trophyp}. 

Further, if $V_1,\ldots,V_k$ are generic hypersurfaces with Newton polytopes $\Delta_1,\ldots,\Delta_k$, then, by the multiplicativity of affine characteristic classes, we have $$\langle V_1\cap\ldots\cap V_k\rangle=\frac{[\Delta_1\rangle\ldots[\Delta_k]}{(1+[\Delta_1])\ldots(1+[\Delta_k])}. \eqno{(1)}$$
In particular, this gives Khovanskii's formula \cite{khov0} for the Euler characteristic of a nondegenerate complete intersection and assigns geometrical meaning to the right hand side of $(1)$ that appeared in \cite{khov0} and subsequent literature as a formal expression.
\end{exa}

\subsection{Polynomial functions of polytopes} \label{Spfp}

This section prepares polyhedral-geometric tools to prove the existence of tropical characteristic classes. Let $P(L)$ be the semigroup of convex polytopes in an $l$-dimensional lattice $L$, and let $[\cdot]:P(L)\to \mathcal{K}^0_{1}(L^*)$ be the inclusion, sending a polytope to its dual fan (see Example \ref{trophyp}). In accordance with this notation, we denote the integer mixed volume of the polytopes $B_1,\ldots,B_n$ in $\Q^n$ by $[B_1]\cdot\ldots\cdot[B_n]$. 

Let 
$e_0,\ldots,e_m$ be the vertices of the standard simplex in $\Q^m$, let $B_0,\ldots,B_m$ and $A_1,\ldots,A_k$ be polytopes in $\Q^n$, and let 
$B$ be the convex hull of the union $\bigcup B_i\times\{e_i\}$ in $\Q^m$.  The volume of the Cayley polytope $B$ can be expressed in terms of mixed volumes of $B_0,\ldots,B_m$ as follows (setting $k=0$):

\begin{lemma} \label{mvcayley} For any polytopes $A_1,\ldots,A_k\subset\Q^m$, we have
$$[A_1]\ldots[A_k]\sum\limits_{b_0+\ldots+b_m=n-k\atop b_0\geqslant 0,\ldots,b_m\geqslant 0}[B_0]^{b_0}\ldots[B_m]^{b_m} = [A_1]\ldots[A_k][B]^{n+m-k}.$$
\end{lemma}
See e.g. Lemma 1.7 in \cite{Edcg} for a proof, based on the Kouchnirenko-Bernstein-Khovanskii formula, or Theorem 24 in \cite{Ejsing} for a combinatorial proof.

\begin{defin} We say that a function $f:P(L)\to\Q$ is a polynomial starting from $S\in P(L)$, if there exist tropical fans $F_i\in \mathcal{K}^0_{i}(L^*)$, such that $f(B)=\sum_i F_i[B]^{n-i}$ for every $B$ that contains $S$. \end{defin}  
\begin{rem} The coefficients of the polynomial are uniquely determined by the values of the function $f$. 
\end{rem}

Let $C$ be the standard $m$-dimensional simplex in $\Q^m$. For a polytope $B\subset\Q^n$, denote the convex hull of the union of $B\times\{0\}$ and $\{0\}\times C$ by $C(B)\subset\Q^n\times\Q^m$. Let $\Sigma\subset(\Q^m)^*\subset(\Q^n\times\Q^m)^*$ be the set of covectors with at least one positive coordinate, i. e. all $\gamma$ such that the support face $C^\gamma$ is disjoint from 0. Recall that the support face $\Delta^\gamma$ of a polytope $\Delta$ with respect to a linear function $\gamma$ is the set of points where the restriction $\gamma:\Delta\to\R$ attain its maximal value. Also recall that faces of two polytopes are said to be compatible, if they support the same covector.
\begin{lemma}\label{startpoly} For every polytope $\Delta\subset\Q^n\times\Q^m$ there exists a polytope $S\subset\Q^n$ such that for every larger $B\supset S$ we have the following: if a face of $C(B)$ disjoint from $B$ is compatible with a face $F\subset\Delta$, then $F=\Delta^\gamma$ for some $\gamma\in\Sigma$.
\end{lemma}
{\sc Proof.} Passing to the projections of the polytopes along $(\delta,\gamma):\Q^n\times\Q^m\to\Q^2$ for finitely many pairs of  $\delta\subset(\Q^m)^*$ and $\gamma\in\Sigma$, we can reduce the question to the case $m=n=1$. This case is obvious: we should take $S$ such that the edges of the triangle $C(S)$ are parallel to the edges of $\Delta$ adjacent to the edge $\Delta^{(0,1)}$. $\quad\square$
  
Let $A\in \mathcal{K}^0_{n-k}(\Z^n\times\Z^m$) be a tropical fan, $k\geqslant 0$.
\begin{utver} \label{poly1} The number $A[C(B)]^{m+k},\, k\geqslant 0,$ is a polynomial $P_A(B)=\sum_i P_{A,i}[B]^{n-i},\, P_{A,i}\in \mathcal{K}^0_i(\Z^n)$, starting from some $S$. This polynomial is multiplicative in $A$: if $A\in \mathcal{K}^{0}(\Z^n\times\Z^m)$ and $A'\in \mathcal{K}^{0}(\Z^{n'}\times\Z^{m'})$, then $P_{A\times A',i}=\sum_j P_{A,j}\times P_{A,i-j}$.
\end{utver}
{\sc Proof.} Since the ring $\mathcal{K}^0$ is isomorphic to the polytope algebra, it is generated by its order 1 component $\mathcal{K}^0_1$, i. e. every fan $A$ can be represented as a linear combination of complete intersections of the form $[A_1]\ldots[A_{n-k}]$ for polytopes $A_i\subset\Q^n\times\Q^m$, so we assume with no loss in generality that $A=[A_1]\ldots[A_{n-k}]$. Let $T_i$ be the projection of $A_i$ to $\Q^m$. For every face $F$ of $T=\sum_i T_i$, represent $F$ as $\sum_i\widetilde F_i$ of faces $\widetilde F_i\subset T_i$, and denote the preimage of $\widetilde F_i$ under the projection $A_i\to T_i$ by $F_i$.

Let $C_I$ be the standard $|I|$-dimensional simplex in $\Q^I\subset\Q^m$, where $\Q^I$ is the coordinate plane, defined by vanishing of the $i$-th coordinates for $i\notin I\subset\{1,\ldots,m\}$. Let $C_I^0$ be its facet, disjoint from $0$, and let $C_I(B)\subset C(B)$ be the convex hull of $B\cup C_I$. For every face $F\subset T$, let $I_F$ be the maximal $I$ such that the faces $C_I^0\subset C$ and $F\subset T$ are compatible (i.e. $F+C_I^0$ is a face of $T+C$). 

Choose any $i_1<\ldots<i_p$ and $j$, then the interiors of the sums $A^{\bar \imath, j}_F=F_{i_1}+\ldots+F_{i_p}+jC_{I_F}(B)$ for any two faces $F\subset T$ do not intersect.
The key observation is as follows: there exists $S\subset\Q^n$ (provided by Lemma \ref{startpoly}), such that for $B\supset S$ the union of $A^{\bar \imath, j}_F$ over all faces $F\subset T$ equals $A_{i_1}+\ldots+A_{i_p}+jC(B)$, so the volume of $A_{i_1}+\ldots+A_{i_p}+jC(B)$ equals the sum of the volumes of $A^{\bar \imath, j}_F$. Combining this observation with the formula $[B_1]\ldots[B_q]=\sum_{i_1<\ldots<i_p}(-1)^p\Vol(B_{i_1}+\ldots+B_{i_p})$ for any polytopes $B_1,\ldots,B_k\subset\R^q$, we conclude that $[A_1]\ldots[A_{n-k}][C(B)]^{m+k}$ equals the sum of $[F_1]\ldots[F_{n-k}][C_{I_F}(B)]^{m+k}$ over all faces $F\subset T$ for $B\supset S$.

Each of the summands $[F_1]\ldots[F_{n-k}][C_{I_F}(B)]^{m+k}$ is a polynomial of $B$, because we can compute it by Lemma \ref{mvcayley} with $B_0=B$, $B_1=C^0_{I_F}$ and $B_{0,1}=C_{I_F}(B)$. This gives an explicit polynomial formula $\sum_i P_{A_1,\ldots,A_{n-k},i}[B]^i$ for $[A_1]\ldots[A_{n-k}][C(B)]^{m+k}$. Also by Lemma \ref{mvcayley}, one can check that $P_{A_1,\ldots,A_{n-k},A'_1,\ldots,A'_{n'-k'},i}=\sum_j P_{A_1,\ldots,A_{n-k},j}\times P_{A'_1,\ldots,A'_{n'-k'},i-j}$ for $A_1,\ldots,A_{n-k}\subset\Q^n\times\Q^m$ and $A'_1,\ldots,A'_{n'-k'}\subset\Q^{n'}\times\Q^{m'}$. $\quad\square$

\subsection{Base points at infinity} \label{Sbasep}

In this section, we discuss families of subvarieties $B_s\subset\CC^n$ that are flat at infinity in the following sense: for every variety $V$ of a complementary codimension, almost all $B_s$ ``do not intersect $V$ at infinity'' (i.e. $V$ and $B_s$ satisfy the assumption of Proposition \ref{tropint0}).

Let $A$ and $B$ be irreducible algebraic varieties of the dimension $a$ and $b$, and $X\subset A\times B$ be a codimension $k$ irreducible subvariety, whose projections $\pi_A:X\to A$ and to $\pi_B:X\to B$ are surjective. For every $s\in A$, denote $\pi_B\pi_A^{-1}(s)$ by $B_s$.
\begin{defin}
The variety $X$ is called the family of subvarieties $B_s\subset B$ (members of the family), parameterised by $s\in A$. A point $y\in B$ is said to be a base point of the family, if $\dim A_y>a-k$.
\end{defin}
The following fact is obvious.
\begin{lemma} \label{nobase0} If the family $X\subset A\times B$ has no base points, then every subvariety $V\subset B$, whose dimension is less than $k$, has empty intersection with a generic member of the family.
\end{lemma}
We are interested in the following special case: $A=\CC^a$ and $B=\CC^b$. 
\begin{defin}\label{defbpinf} The family $X\subset A\times B$ is said to have no base points at infinity, if, for any non-zero $\gamma\in \{0\}\times\Q^b\subset\Q^a\times\Q^b$, the family $\In_{\gamma}X\subset A\times B$ has no base points (see the paragraph before Proposition 
\ref{tropint0} for the meaning of $\In_\gamma$).
\end{defin}
Denote the intersection product of the tropical fan of $X$ and $\{0\}\times\Q^b\subset\Q^a\times\Q^b$ by $L_X$.
\begin{utver} \label{bpinf} 1) The tropical fan of a generic member of the family $X\subset A\times B$ equals $L_X$.

2) If $X\subset A\times B$ has no base points at infinity, and $V\subset B$ is a $k$-dimensional subvariety, then the intersection number of $V$ with a generic member of the family is well defined and equals the tropical intersection number of the tropical fan of $V$ and $L_X$.
\end{utver}
{\sc Proof.} Part (1) is the definition of the product in the ring of conditions 
for $X$ and a fiber of $\pi_A$. If $X$ has no base points at infinity, then, for generic $s\in A$ and any non-zero $\gamma\in\Q^b$, the set $\In_{\gamma} B_s$ is either empty (for $\gamma\notin L_X$), or does not intersect $\In_\gamma V$ by Lemma \ref{nobase0} (for $\gamma\in L_X$, because $X$ has no base points at infinity). Thus, Proposition \ref{tropint0} is applicable to $B_s$ and $V$. 
$\quad\square$

We are especially interested in the following two special cases. Recall that, for a finite subset $A$ in the character lattice of the complex torus $\CC^n$, we denote the space of linear combinations of the characters from $A$ by $\C^A$.
\begin{exa}\label{motiv0} Obviously, the family of planes $\{f\in\C^A\, |\, f(g)=0\}$ for all $g\in\CC^n$ has no base points at infinity.
\end{exa}
\begin{defin} The $a$-logarithmic differential of $f\in\C^A,\, A\subset\Z^n,\, a\in\R^n$, is the map ${\rm dlog}_af=(x_1\frac{\partial f}{\partial x_1}-a_1f,\ldots,x_n\frac{\partial f}{\partial x_n}-a_nf)$.
\end{defin}
\begin{lemma} \label{dlogbp} The graphs of ${\rm dlog}_af$ for $a\in\C^n\setminus($the union of affine spans of the facets of $A)$ and $f\in\C^A$ form a family in $\CC^n\times\CC^n$ with no base points at infinity. 
\end{lemma}
{\sc Proof.} For any $(x,y)=(x_1,\ldots,x_n,y_1,\ldots,y_n)$ and any non-trivial face $F\subset A$, the dimension of the set of all pairs $\Bigl(a=(a_1,\ldots,a_n),\, f=\sum_{b\in A}c_bx^b\Bigr)$, satisfying the equalities $y_i=\sum_{b\in F}(b_i-a_i)c_bx^b,\, i=1,\ldots,n$, does not depend on $(x,y)$.$\quad\square$

\subsection{Existence of characteristic classes} \label{Sconstr1}

We refer to e.g. \cite{khovtor} and \cite{danil} for the background on toric compactifications. For an arbitrary bounded subset $B\subset\R^n$, the space of Laurent polynomials of the form $\sum_{b\in B} c_bx^b$ will be denoted by $\C^B$ instead of $\C^{B\cap\Z^n}$ by a small abuse of notation.

Let $V\subset\CC^n$ be a smooth subvariety.

\begin{defin} The {\bf logarithmic conormal variety} $\mathcal{V}\subset \CC^n\times\CC^n\subset T^*\CC^n$ of $V$ is the set of all $(x,l)$ such that $x\in V$, and the restriction of $\sum l_i\frac{dx_i}{x_i}$ to $T_xV$ equals 0.
\end{defin}

Our idea is to prove that the coefficients of the polynomial $P_{\Trop\mathcal{V}}$ (see Proposition \ref{poly1} for this notation) satisfy the properties of tropical characteristic classes of $V$. The proof starts from the observation that Proposition \ref{bpinf} and Lemma \ref{dlogbp} imply the following.

\begin{lemma} \label{lchar0}  
 
The intersection number of $\mathcal{V}$ and the graph of ${\rm dlog}_a f$ in $\CC^n\times\CC^n$ equals the tropical intersection number of $\Trop\mathcal{V}$ and $[C(B)]^n$.
\end{lemma}

\begin{defin} Denote the Euler characteristic $e(V\cap\{f_1=\ldots=f_i=0\})$ for generic $f_1\in\C^{B_1},\ldots,f_i\in\C^{B_i}$ by $e(V,B_1,\ldots,B_i)$. 
\end{defin}

A sufficient explicit genericity condition for $f_1,\ldots,f_i$ is as follows (the proof is standard).

\begin{lemma} \label{genpos0} Choose a smooth toric compactification of $\CC^n$, whose fan is compatible with the lattice polytopes $B_1,\ldots,B_k$ in $\Q^n$, so that the closure $\overline{M}$ of $M$ intersects the orbits properly (see e.g. \cite{Kaz99} or \cite{Tev}). Let $\mathbf{M}$ be a Whitney stratification of $\overline{M}$, compatible with the orbit stratification of the ambient toric variety. The tuples $(f_1,\ldots,f_k)\in\C^{B_1}\oplus\ldots\oplus\C^{B_k}$, such that the closure of $f_{i_1}=\ldots=f_{i_p}=0$ is smooth and transversal to all strata of $\mathbf{M}$ for every sequence $i_1<\ldots<i_p$, is a non-empty a Zariski open subset in $\C^{B_1}\oplus\ldots\oplus\C^{B_k}$. The Euler characteristic
$e(M\cap\{f_1=\ldots=f_k=0\})$ is the same for all tuples $(f_1,\ldots,f_k)$ in this Zariski open set.
\end{lemma}

The quantities $e(V,B_1,\ldots,B_i)$ satisfy the following higher additivity property (cf. higher additivity for Newton polytopes of discriminants in \cite{Edcg}).

\begin{lemma} \label{hiadd1} We have $e(M,B_0+B'_0,B_1,\ldots,B_k)+2e(M,B_0,B'_0,B_1,\ldots,B_k)=e(M,B_0,B_1,\ldots,B_k)+e(M,B'_0,B_1,\ldots,B_k)+e(M,B_0+B'_0,B_0,B'_0,B_1,\ldots,B_k)$.
\end{lemma}
{\sc Proof.} Compare the Euler characteristic of $M\cap\{gh=f_1=\ldots=f_k=0\}$ and its smoothening $M\cap\{gh+\varepsilon f=f_1=\ldots=f_k=0\}$ for generic $g\in\C^{B_0},\, h\in\C^{B'_0},\, f\in\C^{B_0+B'_0}$. $\quad\square$

\begin{utver} \label{char1} If $\mathcal{V}$ is not empty, then there exists a polytope $S\subset\Q^n$, such that, for any $n$-dimensional polytope $B\supset S$, we have $$e(V,B)=e(V)-(-1)^{\dim V}P_{\Trop\mathcal{V}}(B)=e(V)-(-1)^{\dim V}\sum_{j}P_{\Trop\mathcal{V},j}[B]^{n-j}.$$
\end{utver}
{\sc Proof.} For a generic $f\in\C^B$ in the sense of Lemma \ref{genpos0} and $a$ in the interior of $B$, by the standard Morse theory argument for the function $h:V\to\R,\, h(x)=|f(x)/x^a|$, we have $e(V,B)=e(V)-(-1)^{\dim V}\cdot($the number of critical points of $h)$. The key observation is that the number of critical points equals the intersection number of $\mathcal{V}$ and the graph of ${\rm dlog}_a f$. The latter, by Lemma \ref{lchar0}, equals the tropical intersection number of $\Trop\mathcal{V}$ and $[C(B)]^n$, which equals $P_{\Trop\mathcal{V}}(B)$ starting from some polytope by Proposition \ref{poly1}. $\quad\square$ 

We can drop the assumption of non-emptiness of the conormal variety as follows.

\begin{utver} \label{char15}
1) If $\mathcal{V}$ is empty, then $V$ is preserved by the action of a certain non-trivial subtorus of $\CC^n$.

2) Let $\CC^k$ be the maximal subtorus of $\CC^n$, preserving $V$, so that $V=V'\times\CC^k$, then there exists a polytope $S\subset\Q^{n-k}$, such that we have $e(V,B)=-(-1)^{\dim V+k}P_{\Trop\mathcal{V'}}(B)$ for any $n$-dimensional $B\supset S$.
 
\end{utver}
The proof is the same as for Proposition \ref{char1}: let $x_1,\ldots,x_n$ be monomial coordinates on $\CC^n$, such that $\CC^k=\{x_{k+1}=\ldots=x_n=1\}$, then the set 
of critical points of $h$ can be represented as the intersection of $\mathcal{V'}\times\CC^k$ and the set $\{(x,y)\, |\, y_i=x_i\frac{\partial f}{\partial x_i}-a_if$ for $i\leqslant k$ and $\frac{\partial f}{\partial x_i}=0$ for $i>k\}$.

\begin{utver} \label{char2} There exists a polytope $S\subset\Q^n$, such that, for 
polytopes $B_0,\ldots,B_m\supset S$, we have $e(V,B_0,\ldots,B_m)=$

\flushright $e(V)-(-1)^{\dim V}\Bigl(P_{\Trop\mathcal{V},n}+(-1)^m\sum_{j<n-m}P_{\Trop\mathcal{V},j}\sum\limits_{b_0+\ldots+b_m=n-j\atop b_0\in\N,\ldots,b_m\in\N}[B_0]^{b_0}\ldots[B_m]^{b_m}\Bigr)$.
\end{utver}
{\sc Proof.} In the notation of Lemma \ref{mvcayley}, let $f_I$ be the polynomial $\sum_{i\in I} \lambda_if_i$ on $\CC^n\times\CC^I$, and let $B_I$ be its Newton polytope. Consider the action of the complex circle $\CC$ on $\CC^I$ by coordinatewise multiplication and the projection $\pi_I:\CC^n\times\CC^I\to\CC^n$. By the additivity of the Euler characteristic, the desired $e(V\cap\{f_0=\ldots=f_m=0\})$ equals $\sum_I e\Bigl(\{f_I=0\}\cap (V\times\CC^n)/\CC\Bigr)-me(V)$, because the restriction of the projection $\CP^m\times V\to V$ to $\{\sum_i\lambda_if_i\}/\CC$ has fibers $\CP^{m-1}$ over the points of $\{f_0=\ldots=f_m=0\}$ and $\CP^m$ over the other points of $V$.

We can now compute $ e\Bigl(\{f_I=0\}\cap (V\times\CC^n)/\CC\Bigr)$ for every $I$ by Proposition \ref{char1} or \ref{char15} in terms of $P_{\Trop\mathcal{V}}(B_I)$, and simplify the answer by Lemma \ref{mvcayley}. $\quad\square$

\begin{utver} \label{char0n} If $\mathcal{V}$ is non-empty (i. e. $V\subset\CC^n$ is not preserved by non-trivial tori), then the leading term of the polynomial $P_{\Trop\mathcal{V}}(B)$ equals $\Trop(V)[B]^{\dim V}$, and the constant term equals $(-1)^{\dim V}e(V)$. 
\end{utver}
{\sc Proof.} The equalities $P_{\Trop\mathcal{V},j}=0,\, j<n-\dim V$, $P_{\Trop\mathcal{V},n-\dim V}=\Trop V$ and $P_{\Trop\mathcal{V},n}=(-1)^{\dim V}e(V)$ follow from Proposition \ref{char2} for $m=\dim V-1$, because $e(V,B_0,\ldots,B_m)$ equals $[B_0]\ldots[B_m]\Trop V$ for all $B_0,\ldots,B_m$. $\quad\square$

\begin{utver} \label{char25} The equality of Proposition \ref{char2} is valid without the assumptions $\dim B=n$ and $B\supset S$.
\end{utver}
{\sc Proof.} For $m=\dim V-1$, the statement follows from Proposition \ref{char0n}. Assume that we have proved the statement for all $m$ greater than a given value $m_0$. Then we can prove it for $m=m_0$ by induction on the number $p$ of the polytopes $B_i$ that do not contain $S$ or have dimension less than $n$. For $p=0$, the desired statement coincides with the one of Proposition \ref{char2}. For $p>0$, apply Lemma \ref{hiadd1} with $k=m$, $B_{0},\ldots,B_{p-1}$ containing $S$ and $n$-dimensional, and $B'_0,B_{p},\ldots,B_{m}$ arbitrary. Note that $B_0+B'_0$ is also $n$-dimensional and containing $S$, so all the terms in the equality of Lemma \ref{hiadd1}, except for $e(M,B'_0,B_1,\ldots,B_m)$, can be computed by the inductive assumption. This can be used to compute $e(M,B'_0,B_1,\ldots,B_m)$ and leads to the desired formula for it. $\quad\square$

\begin{defin} \label{finchar} Let $V:\CC^n\to\R$ be a constructible function (that is, a linear combination of characteristic functions of algebraic subsets). Define its characteristic class $\langle V\rangle=\sum_i\langle V\rangle_i,\, \langle V\rangle_i\in \mathcal{K}^0_i(\Z^k)$, as follows.

1) If $V$ is (the characteristic function of) a smooth subvariety $W\subset\CC^n$, then define $\langle W\rangle_i$ as $(-1)^{i+n+\dim W}P_{\Trop\mathcal{W},i}$.

2) If $V$ is (the characteristic function of) a smooth subset of $\{f\ne 0\},\, f:\CC^n\to\C$, 
then define $\langle V\rangle_i$ as the image of $\langle V\times\CC\cap($the graph of $f)\rangle_i$ under the projection $\CC^n\times\CC\to\CC^n$.

3) For arbitrary $V$, represent it as a linear combination $\sum_i\alpha_jV_j$ of functions $V_i$ of the form (2), and define $\langle V\rangle_i$ as $\sum_j\alpha_j\langle V_i\rangle_j$.
\end{defin}

\begin{theor} \label{fincharth} 1) For all lattice polytopes $B_1,\ldots,B_k$ in $\Q^n$ and generic $f_i\in\C^{B_i}$, the Euler characteristic integral of $V$ over $\{f_1=\ldots=f_k=0\}$ (see e.g. \cite{mcph} and \cite{viro}) equals $(-1)^{\dim V+k}\sum_{j\leqslant n-k}\langle V\rangle_j\sum_{b_1+\ldots+b_k=n-j\atop b_1\in\N,\ldots,b_k\in\N}[B_1]^{b_1}\ldots[B_k]^{b_k}$.

2) We have $\langle U\times V\rangle=\langle U\rangle\times\langle V\rangle$.

3) $\langle V\rangle_n$ equals the Euler characteristic integral of $V$.

4) $\langle V\rangle$ does not depend on the decomposition of $V$ into a linear combination of characteristic functions, chosen in Definition \ref{finchar}.

5) For generic $c\in\CC^n$, we have $\langle U\cap(cV)\rangle=\langle U\rangle\langle V\rangle$.

6) For a generic hypersurface $V$ with Newton polytope $\Delta$, we have $\langle V\rangle=\frac{[\Delta]}{1+[\Delta]}$.
\end{theor}
{\sc Proof.} For characteristic functions of smooth subvarieties of tori, (1) is Proposition \ref{char25}, (3) is Proposition \ref{char0n}, and (2) follows from multiplicativity of $P_A(B)$ (see Proposition \ref{poly1}). For arbitrary characteristic functions, (1,2,3) follow by the additivity of Euler characteristic. They imply (4), because $\langle V\rangle_i,\, i<n$, are uniquely determined by (1) for $k=1$, and $\langle V\rangle_n$ is uniquely determined by (3). To prove (5), represent $U\cap(cV)$ as $U\times V\cap\{(x,y)\, |\, x_i-c_iy_i=0\}$, then the desired equality follows from (1,2). Part (6) follows from Proposition \ref{char1}, because, computing $e(V,B)$ in its statement by the Kouchnirenko-Bernstein-Khovanskii formula \cite{khov0}, we have $P_{\Trop\mathcal{V},j}=[\Delta]^j$. $\quad\square$

We now prove the relation of tropical characteristic classes to Schwarz-MacPherson classes (Part 6 of Definition \ref{defspher} and Proposition \ref{schoncc}). For this, recall the 
description of the cohomology ring of a smooth toric variety in terms of tropical fans. Let $\Sigma$ be a $\Z$-simple fan, and let $X_\Sigma$ be the corresponding smooth toric compactification of the torus $\CC^n$. We say that a $k$-dimensional tropical fan $(P,\varphi)$ is {\bf compatible} with $\Sigma$, if its support set $P$ is contained in the union of the $k$-dimensional cones of $\Sigma$.
\begin{utver}[\cite{brion}, \cite{FS}]\label{homoltor} The cohomology ring of $X_\Sigma$ is naturally isomorphic to the ring of tropical fans compatible with $\Sigma$. 
The isomorphism sends a codimension $k$ cohomology class $\alpha$ to the tropical fan $(P,\varphi)$, where $P$ is the union of the $k$-dimensional cones of $\Sigma$, and the value of $\varphi$ on every such cone $C$ equals the restriction of $\alpha$ to the closure of the $C$-orbit in $X_\Sigma$.
\end{utver}

We identify the cohomology and homology of $X_\Sigma$ with the Poincare duality. For a subvariety $V\subset\CC^n$, we consider the Schwarz-MacPherson class $c^{SM}_{i,\Sigma}(V)\subset H^i(X_\Sigma)$ of a non-closed set $V\subset X_\Sigma$ as a tropical fan, compatible with $\Sigma$. 
\begin{theor} \label{thm239} 1) If the closure of $V$ in $X_\Sigma$ is smooth and transversal to the orbits, then the multiplicity of every cone $C\in\Sigma$ in the tropical fan $c^{SM}_{i,\Sigma}(V)$ equals the Euler characteristic of {\rm (closure of $V)\cap(C$-orbit)}.

2) If $U\subset\CC^n$ is arbitrary, $V\subset\CC^n$ is a generic complete intersection, defined by the equations with Newton polytopes $N_1,\ldots,N_k$, and the dual fans $[N_1],\ldots,[N_k]$ are compatible with $\Sigma$, then $c^{SM}_\Sigma(U\cap V)=c^{SM}_\Sigma(U)c^{SM}_\Sigma(V)$.

3) For arbitrary $U\subset\CC^n$, choose $\Sigma$ compatible with the tropical characteristic class $\langle U\rangle$, constructed above. 
Then $\langle U\rangle$ equals $c^{SM}_\Sigma(U)$. 
\end{theor}
\begin{lemma}\label{transvsm} (see Theorem 3.1 in \cite{aluf})
Let $M$ be a smooth compact variety, $V\subset M$ be an arbitrary subvariety, and $D$ be a smooth hypersurface, transversal to $V$. Then
1) $c^{SM}(V\cap D)=c^{SM}(V)\cdot c^{SM}(D)$ in the cohomology ring $H^\bullet(M)$,

2) $c^{SM}(V\cap D)\in H^\bullet(D)$ is equal to the restriction of $c^{SM}(V\setminus D)\in H^\bullet(M)$.
\end{lemma}
{\it Proof of Theorem \ref{thm239}.} Part 1 follows from Lemma \ref{transvsm}(2) and Proposition \ref{homoltor}, because the highest Schwarz-MacPherson class equals the Euler characteristic. Part 2 follows from Lemma \ref{transvsm}(1) for a generic hypersurface $V$ and by induction for a generic complete intersection $V$. 

Part 3 for a generic complete intersection follows from Part 1 and Theorem \ref{fincharth}(6). For an arbitrary $U$, choose a generic complete intersection $V$, such that the dual fans of the Newton polytopes of its equations are compatible with $\Sigma$. Then we have $c^{SM}_\Sigma(V)=\langle V\rangle$, and the highest components of $c^{SM}_\Sigma(U)\cdot c^{SM}_\Sigma(V)$ and $\langle U\rangle\cdot\langle V\rangle$ are both equal to the Euler characteristic of $U\cap V$ (by Part 2 and Theorem \ref{fincharth}(5) respectively). This equality for arbitrary $V$ implies that $\langle U\rangle=c^{SM}_\Sigma(U)$, because the cohomology ring of $X_\Sigma$ is generated by complete intersections, and the multiplication of cohomologies of complementary dimensions is a nondegenerate bilinear form. $\quad\Box$

\subsection{Characteristic classes of critical complete intersections} \label{Scci}

Non-degenerate complete intersections are the simplest example of sch\"{o}n varieties. We shall need the next simplest example: the set of critical points of a projection $\{f=0\}\subset\CC^n\stackrel{p}{\to}\CC^k$, where $f$ is a Laurent polynomial with generic coefficients, and $p$ is a morphism of complex tori. This complete intersection is degenerate, but still sch\"{o}n, as we prove below.

Denote the $(n+1)$-dimensional space of affine linear functions on $\Z^n$ by $(\Z^n)^\star$, to be distinguished from the dual space $(\Z^n)^*$. For every subspace $L\subset\Z^n$, let $L^\perp\subset (\Z^n)^\star$ be the set of affine functions on $\Z^n$ that vanish on $L$. For a Laurent polynomial $f(x)=\sum_{b\in\Z^n} c_bx^b$ on $\CC^n$ with the standard coordinates $x=(x_1,\ldots,x_n)$ and an element $\alpha\in(\Z^n)^\star$, denote the polynomial $\sum_{b\in\Z^n} \alpha(b)c_bx^b$ by $\partial_\alpha f$. 
\begin{exa} If $\alpha$ is the $i$-th coordinate function, then $\partial_\alpha f=x_i\frac{\partial f}{\partial x_i}$.
\end{exa}
\begin{defin} The set  $\{\partial_D f=0\}\subset\CC^n$, given by the equations $\partial_\alpha f=0,\, \alpha\in D$, is called a {\bf critical complete intersection}.
\end{defin}
\begin{utver} \label{critschon} 1) For every finite $B\subset\Z^n$, every vector subspace $D\subset (\Z^n)^\star$ and a generic Laurent polynomial $f(x)=\sum_{b\in B} c_bx^b$,
the critical complete intersection $\{\partial_D f=0\}$ is indeed a regular complete intersection.

2) Moreover, it is sch\"{o}n.
\end{utver}
Note that the sch\"{o}n complete intersection $\{\partial_D f=0\}$ is not in general non-degenerate with respect to its Newton polytope (\cite{khovtor}) even for generic $f$: for instance, if $f\in\C[x,y]$ is a generic polynomial with a given Newton polygon $N$, then the restrictions of the equations $f=y\frac{\partial f}{\partial y}=0$ to each of the horizontal edges of $N$ are compatible, because the restriction of $y\frac{\partial f}{\partial y}$ is proportional to the restriction of $f$.

Part 1 of Proposition \ref{critschon} is a standard corollary of the Bertini---Sard theorem, see the subsequent paragraph. The proof of Part 2, including an explicit construction of a sch\"{o}n compactification, occupies the rest of this subsection.

{\sc Proof of Part 1.} With no loss in generality we can assume that $B$ is not contained in an affine hyperplane. Then we can choose a basis $\alpha_0,\ldots,\alpha_k$ in $D$ and points $b_0,\ldots,b_k$ in $B$ such that the matrix $\alpha_i(b_j)$ is triangular and non-degenerate, i.e. $\partial_{\alpha_i} f$ does not depend on $c_{b_{i+1}},\ldots,c_{b_{k}}$. Then, by induction on $i$, we can prove that $\partial_{\alpha_0} f=\ldots=\partial_{\alpha_i} f=0$ is a regular complete intersection for almost all $f$: for every $i$, take $c_{b_i}$ to be a regular value of the restriction of the function 
$-\sum_{b\in B\setminus\{b_i\}} \frac{\alpha_i(b)}{\alpha_i(b_i)}c_bx^{b-b_i}$ to the regular set $\partial_{\alpha_0} f=\ldots=\partial_{\alpha_{i-1}} f=0$.
$\quad\Box$ 

The proof of Part 2 is based on an explicit construction of a sch\"{o}n compactification for $\{\partial_D f=0\}$. We shall use the following notation to describe it. 
For every $\gamma\in(\Z^n)^\star$ and $B\subset\Z^n$, let $B^\gamma$ be the set of points at which the restriction of $\gamma$ to $B$ attains its maximum. For every sequence of sets $B_k\subset\Z^n$, denote the {affine span} of $\bigcup_{i<k}B_k$ by $B_{<k}$.

\begin{defin} 1) For every $\gamma\in(\Z^n)^\star$, define a finite sequence of subsets $B^\gamma_k\subset B$ inductively by $k\geqslant 0$ as follows: the set $B^\gamma_k$ consists of all points at which the restriction of the function $\gamma$ to $B\setminus B^\gamma_{<k}$ attains its maximal value, unless $(B^\gamma_{<k})^\perp\cap D=0$,
and $B_k$ is not defined in the latter case.

2) For a finite sequence of non-empty non-intersecting subsets $B_k\subset B$, the cone $C_{B_\bullet}$ is the set of all $\gamma\in (\Z^n)^*$ such that the sequence $B^{\gamma}_k$ equals the sequence $B_k$. The fan that consists of the cones $C_{B_\bullet}$, as $B_\bullet$ runs over all sequences of 
subsets of $B$, will be denoted by $C_B$.
\end{defin}

\begin{rem} The codimension of a non-empty cone $C_{B_\bullet}$ equals the dimension of $\sum_i B_i$.
\end{rem}

\begin{exa} If $B=\{0, e_1, e_2, e_3\}$ is the standard simplex in $\R^3$, and $D$ is generated by the functions $1$ and $x+y$, then $C_B$ is the minimal subdivision of the dual fan of $B$ by the line generated by $x+y$. For example, the 1-dimensional cone generated by $-x-y$ corresponds to the sequence $B_{\bullet}=(B_0, B_1)=(\{0, e_1\},\, \{e_2, e_3\})$.
\end{exa}

Now Proposition \ref{critschon}(2) can be restated in more detail as follows. For every sequence of non-intersecting sets $B_k\subset\Z^n$, choose a direct sum decomposition $D=\bigoplus_i B^i$ such that $\bigoplus_{i\geqslant k} B^i$ equals $B_{<k}^\perp\cap D$ for every $k$. For every subset $B'\subset B$ and a polynomial $f(x)=\sum_{b\in B}c_bx^b$, denote the {\bf restriction of $f$ to $B'$} by $f^{B'}(x)=\sum_{b\in B'}c_bx^b$.

Choose any subdivision $\widetilde C_B$ of the fan $C_B$ such that the $\widetilde C_B$-toric variety is smooth. For every cone $C_{B_\bullet}$, the orbits of  the $\widetilde C_B$-toric variety, corresponding to subcones of $C_{B_\bullet}$, will be called the $C_{B_\bullet}$-orbits.
\begin{utver} \label{compcci} 1) For almost all $f$ of the form $\sum_{b\in B}c_bx^b$, the closure of the critical complete intersection $\{\partial_D f=0\}$ in the $\widetilde C_B$-toric variety is smooth and intersects its orbits transversally. 

2) The closure of $\{\partial_D f=0\}$ intersects a $C_{B_\bullet}$-orbit by the same set as the closure of the set given by the equations $\partial_{B^i}f^{B_i}=0$ for all $i$.
\end{utver}
{\sc Proof.} 
Choose a basis $\alpha_{0i},\ldots,\alpha_{k_ii}$ in $B_i$. The union of these bases is a basis in $D$. For every linear function $l\in C_{B_\bullet}$, the maximum of $l$ on the Newton polytope of the polynomial $\partial_{\alpha_{ji}}f$ is attained on its face $B_i$. Thus, the intersection of a $C_{B_\bullet}$-orbit and the closure of $\{\partial_D f=0\}$ is contained in the intersection of the $C_{B_\bullet}$-orbit and the closure of $\cap_{i,j}\{\partial_{\alpha_{ji}} f=0\}$, and contains the regular part of the latter intersection.
However, the latter intersection is regular for generic $f$, because, for every $i$, the equations $\partial_{B^i}f^{B_i}=0$ define a regular set $S_i$ by Proposition \ref{critschon}(1), and, for generic $c_i\in\CC^n$, the sets $c_i S_i$ intersect transversally. Thus we have proved that the intersection of the closure of $\{\partial_D f=0\}$ with every $C_{B_\bullet}$-orbit is regular for every cone $C_{B_\bullet}$, and is given by the equations $\partial_{B^i}f^{B_i}=0$. This concludes the proof of Propositions \ref{compcci} and \ref{critschon}.
$\quad\Box$

Proposition \ref{schoncc} reduces computation of tropical characteristic classes of critical complete intersections to computation of the Euler characteristics of critical complete intersections. For every sequence $B_\bullet$, let $Z_{B_\bullet}$ be $\Bigl($the set given by the equations $\partial_{B^i}f^{B_i}=0$ for all $i\Bigr)\Bigl/\Bigl($the subtorus of $\CC^n$ generated by $t^\gamma$ over all $\gamma\in C_{B_\bullet}\Bigr)$.

\begin{sledst} The $j$-th tropical characteristic class of $\{\partial_D f=0\}$ consists of the cones $C_{B_\bullet}$ over all sequences $B_\bullet$ such that $j=\dim\sum_i B_i$, with the multiplicities equal to the Euler characteristics of $Z_{B_\bullet}$.
\end{sledst}

\begin{rem} For most of sequences $B_\bullet$, the cone $C_{B_\bullet}$ is either empty or has multiplicity 0. For example, $C_{B_\bullet}$ is empty unless the sets $B_i$ are non-empty and non-overlapping 
and $B_0$ is contained in the boundary of the convex hull of $B$. The multiplicity of $C_{B_\bullet}$ also equals 0 if $\dim D\cap B_{<k}^\perp>\dim\sum_{i<k}B_i$ for some $k$.
\end{rem}

The sets of the form $Z_{B_\bullet}$ that we face in the framework of this paper are very simple (points and lines). However we notice that such Euler characteristics can be computed in the general case by methods of \cite{Edcg}. For example, we give the answer for 0--dimensional critical complete intersections.

{\sc Notation.} Assume that $D\subset(\Z^n)^\star$ is a hyperplane. It is uniquely determined by the common zero $z\in\RP^n$ of all affine linear functions $\alpha\in D$. Let $\mathfrak{A}$ be the set $\{\Gamma\,|\,\Gamma$ is a face of the convex hull of $A$, and the affine span of $\Gamma$ contains $z\}$. Let $N(A,D)$ be the sum $\sum_{\Gamma\in\mathfrak{A}} e^\Gamma_A \Vol\Gamma$, where $e^\Gamma_A$ is the Euler obstruction of $A$ at $\Gamma$, see \cite{Edcg} for the definition.
\begin{lemma} \label{dim0crit} If $A$ generates $\Z^n$ and $D$ is a hyperplane, then, for a generic Laurent polynomial $ f(x)=\sum_{a\in A} c_ax^a$, the 0-dimensional critical complete intersection $\{\partial_D f=0\}$ consists of $N(A,D)$ points.
\end{lemma}
This is the Bernstein formula \cite{bernst0} with some roots hidden at infinity (more precisely, at the orbits of the $A$-toric variety, corresponding to the faces from $\mathfrak{A}$), see \cite{Edcg} for details.

\section{Affine Pl\"ucker formulas} \label{SSapf}

The problematics of affine multisingularity theory is introduced in Section \ref{Stmt}. Its relation to tropical enumerative geometry is outlined in Section \ref{Senum}. The affine version of the Pl\"ucker formulas is obtained in Sections \ref{Sdiscr}---\ref{Smain}: the tropical fan of the universal singularity stratum of  codimension 1 (i.e. the Newton polytope of the $A$-discriminant for a non-dual defective set $A\subset\Z^n$) is computed in Section \ref{Sdiscr}; additional assumptions on $A$ ($k$-versality, Definition \ref{defvers}), playing the role of dual defectiveness for the codimension 2 universal multisingularity strata, are introduced in Section \ref{Smain0}; the tropical fans of these strata are described in Section \ref{Smain} (see Section \ref{Sproofs123} for a proof of this description). Section \ref{Sred} explains how to specialize these results to various non-universal settings such as generic polynomial maps of generic affine hypersurfaces, and an example of the computation for a simple special case is given in Section \ref{Sexample}. Application of affine Pl\"ucker formulas requires somewhat tedious computations with secondary fans that can be simplified with a formula for the volume of a fiber body, presented in Section \ref{Svfp}.

\subsection{Affine multisingularity theory} \label{Stmt}

A singularity theory is an equivalence relation on the set of germs of algebraic varieties. An equivalence class is called a singularity. 
\begin{exa} \label{xxv} We shall study $\mathcal{A}_1$ and $\mathcal{A}_2$ singularities. The $\mathcal{A}_i$ singularity is the class of hypersurface germs in $\C^n$ that can be represented as $x_1^{i+1}+x_2^2+\ldots+x_n^2=0$ in suitable local coordinates. In particular, for $n=2$, the $\mathcal{A}_1$ and $\mathcal{A}_2$ singularities are nodes and cusps of plane curves, and we shall denote them by $\times$ and $\prec$ respectively.
\end{exa}
A finite tuple of isolated singularities $S=(S_1,\ldots,S_k)$ is called a multisingularity. The $S$-multisingularity stratum of a proper morphism of varieties $N\to M$ is the set of all $x\in M$ whose preimages have exactly $k$ singularities $S_1,\ldots,S_k$. One of the first goals of singularity theory is to describe the topology (e. g. the homology class) of a multisingularity stratum in terms of the topology (e. g. the characteristic classes) of $M$ and $N$, provided that the morphism is ``generic enough''. 

However, this theory does not suit the study of multisingularity strata of polynomial mappings $f:\C^n\to\C^m$ or $\CC^n\to\CC^m$, as explained in the introduction. 
A possible way to understand singularity theory in this affine setting is as follows: a generic mapping $f=(f_1,\ldots,f_m)$ is a mapping whose components $f_i$ are generic linear combinations of given monomials, and the goal of our singularity theory is to describe the cohomology classes of multisingularity strata of the map $f$ in the ring of conditions of $\CC^n$, that is, in the ring of tropical fans $\mathcal{K}^0(\Z^n)$. 

Instead of studying individual maps, we first study the {\bf universal setting}: identify elements of $\Z^n$ with monomials in the variables $x_1,\ldots,x_n$, let $\C^A$ be the space of linear combinations of monomials from a finite set $A\subset\Z^n$, choose finite sets $A_1,\ldots,A_m$ in $\Z^n$.
The universal $S$-multisingularity stratum is the set of all tuples of polynomials $f\in\C^{A_1}\oplus\ldots\oplus\C^{A_m}$,
such that the variety $f=0$ has exactly $k$ singularities $S_1,\ldots,S_k$, and our goal is to describe tropical fans of universal multisingularity strata in terms of the sets $A_1,\ldots,A_m$. 
Then we explain how to reduce individual problems regarding generic mappings $f:\C^n\to\C^m$ or $\CC^n\to\CC^m$ to the universal case.

A lot of recent research can be considered from this point of view: Newton polytopes of $A$-discriminants,  sparse resultants and mixed discriminants in \cite{GKZ94}, \cite{S94}, \cite{GP}, \cite{EK}, \cite{STY}, \cite{DFS}, \cite{ST}, \cite{Edcg}, \cite{Eadv}, \cite{Smix}, \cite{Dmix}, \cite{Tsikh} are tropical fans of the $(\times)$-singularity stratum. The tropical fan of higher dual toric varieties \cite{Dhi} is the tropical fan of the stratum of hypersurfaces with a single singularity of a given order. 

If the sets $A_1,\ldots,A_m$ are ``small enough'', then collections of isolated singularities do not appear at all in the fibers $f=0,\; f\in  \C^{A_1}\oplus\ldots\oplus\C^{A_m}$, and one can be interested instead in the strata $\Sigma_k$ of all $f$ such that the singular locus of $f$ has a given dimension $k$. For $m=1$, the lowest codimension non-empty stratum $\Sigma_k$ is the $A$-discriminant variety, and its tropical fan is computed in \cite{DFS} (for any $k$). If $A_1=\ldots=A_m$ is the standard simplex, then the strata $\Sigma_k$ are determinantal varieties. In \cite{Shitov}, all triples $(m,n,k)$ are classified, for which the minors of the matrix form a tropical basis of the ideal of the determinantal variety. This, in particular, gives the tropical fan of $\Sigma_k$. For other triples $(m,n,k)$, the tropical fan is unknown. 

Tropical enumeration of plane curves with nodes and cusps in \cite{M}, \cite{MB}, \cite{Sh}, \cite{Shmult} can also be seen as a partial computation of tropical fans of corresponding Severi varieties, that is, the universal $(\times,\ldots,\times,\prec,\ldots,\prec)$-multisingularity strata, see Section \ref{Senum} and \cite{Kstud} for details.

Summarizing this, the only singularity strata, whose complete tropical fans are known, are

-- the lowest codimension singularity strata,

-- the single singularity strata for singularities of a given order,  

-- the strata with known tropical bases of their defining ideals.

We shall compute the tropical fan of the simplest possible strata outside of this list, namely, the codimension 2 multisingularity strata: the $(\times\times)$-stratum of all $f\in\C^A$ with two nondegenerate singular points and the $(\prec)$-stratum of all $f\in\C^A$ with a single minimally degenerate singular point.

Many relations for tropical characteristic classes of multisingularity strata follow directly from additivity of the Euler characteristic. For example, let $X$ be a smooth toric variety, and let $V\subset X\times\CC^N$ be a subvariety. For every multisingularity $S=(S_1,\ldots,S_k)$, let $e(S)$ be the sum of the reduced Euler characteristics of the Milnor fibers of the isolated singularities $S_1,\ldots,S_k$, let $\langle S\rangle$ be the tropical characteristic class of the $S$-multisingularity stratum of the projection $\pi:V\to\CC^N$, and let $p$ be the codimension of the set of all points $x\in\CC^N$, such that the fiber $\pi^{-1}(x)$ has a non-isolated singularity. If a generic fiber of $\pi$ is smooth, then, denoting its Euler characteristic by $e_0$, we have for $i<p$:
\begin{utver}\label{eqstar}  
$$e_0\cdot\langle \pi(V)\rangle_i-\sum_{S} e(S) \cdot\langle S\rangle_i=\sum_{T\subset X} \pi_*\langle V\cap T\rangle_{i+\dim T}, \eqno{(*)}$$ 
where $T$ runs over all orbits of $X$, and $S$ runs over all (isolated) multisingularities. 
\end{utver}

{\sc Proof.} Choose an arbitrary subvariety $Y\subset\CC^N,\, \dim Y<p$,  and a generic element $g\in\CC^N$. Counting the Euler characteristic of $\pi^{-1}(gY)$ fiberwise, it equals $\sum_S \bigl(e_0-e(S)\bigr)e(gY\cap S)$, because the Euler characteristic of every fiber $\pi^{-1}(y),\, y\in S$, equals $e_0-e(S)$. Since $e\bigl(\pi^{-1}(gY)\bigr)=\sum_{T} e(\pi^{-1}(gY)\cap T)=\sum_{i,T} \langle Y\rangle_{N-i}\pi_*\langle V\cap T\rangle_{i+\dim T}$, and $e(gY\cap S)=\sum_i \langle Y\rangle_{N-i}\langle S\rangle_{i}$, we have $$\sum_i \langle Y\rangle_{N-i} \Bigl( e_0\cdot\langle \pi(V)\rangle_i-\sum_{S} e(S) \cdot\langle S\rangle_i \Bigr) = \sum_i \langle Y\rangle_{N-i} \Bigl( \sum_{T\subset X} \pi_*\langle V\cap T\rangle_{i+\dim T} \Bigr)$$ for arbitrary $Y$. By the nondegeneracy of the intersection number pairing in the ring of tropical fans, this is equivalent to $(*)$. $\quad\Box$

\subsection{Relation to tropical correspondence theorems} \label{Senum}

In this subsection, we restate tropical correspondence theorems in terms of intersection numbers of tropical fans of multisingularity strata. This allows us to consider tropical enumeration as the classical approach to enumerative geometry, up to substitution of conventional cohomology rings, Chern classes and Thom polynomials with the ring of tropical fans, tropical characteristic classes and tropical fans of multisingularity strata.

Recall that we write tropical objects in Fraktur. Every codimension 1 tropical fan with positive multiplicities equals the hypersurface $\mathfrak{f=0}$ for some tropical polynomial $\mathfrak{f}$ with unit coefficients. Conversely, for every tropical hypersurface $\{\mathfrak{x\, |\, f(x)=0}\}$ and every point $\mathfrak{x}_0$ in $\TT^m$, the shifted hypersurface $\{ \mathfrak{x\, |\, f(x/x}_0)=\mathfrak{0}\}$ coincides in a small neighborhood of $(\mathfrak{1,\ldots,1})$ with a certain codimension 1 tropical fan that we denote by $\mathfrak{\{f=0\}}_{\mathfrak{x}_0}\in\mathcal{K}^0_1(\Z^m)$. 
\begin{defin} \label{deftropint1} For any collection of tropical fans and hypersurfaces $\mathfrak{H}_1,\ldots,\mathfrak{H}_k$ in $\TT^m$, whose codimensions sum up to $m$, the number $(\mathfrak{H}_1)_\mathfrak{x}\cdot\ldots\cdot(\mathfrak{H}_k)_\mathfrak{x}\in\mathcal{K}^0_m(\Z^m)=\Q$ is called the {\bf tropical intersection number} of $\mathfrak{H}_1,\ldots,\mathfrak{H}_k$ at the point $\mathfrak{x}\in\TT^m$.
\end{defin}
Recall that the following {\bf displacement rule} allows us to count the tropical intersection number and can be taken as its definition. For a rational affine subspace $L\subset\R^m$, by a slight abuse of notation, denote the lattice $\{l\in(\Z^m)^* \, |\, l(L)={\rm const}\}$ by $L^\perp$.
\begin{utver}
1) Assume that $\mathfrak{H}_1,\ldots,\mathfrak{H}_k$ intersect transversally at $\mathfrak{x}$, i.e. $\mathfrak{H}_i\subset\TT^m=\R^m$ in a small neighborhood of $\mathfrak{x}$ coincides with a rational affine space $L_i\ni \mathfrak{x}$ endowed with a multiplicity $m_i$, and $\sum_i L_i^\perp$ is an $m$-dimensional lattice $\mathcal{L}$. Then the tropical intersection number of $\mathfrak{H}_1,\ldots,\mathfrak{H}_k$ at $\mathfrak{x}$ equals $m_1\ldots m_k \cdot\bigl|(\Z^m)^*/\mathcal{L}\bigr|$.

2) In the general case, let $\mathfrak{\widetilde H}_i$ be the copy of $\mathfrak{H}_i$ shifted by a small generic vector. Then $\mathfrak{\widetilde H}_1,\ldots,\mathfrak{\widetilde H}_k$ intersect at finitely many points $\mathfrak{x}_1,\ldots,\mathfrak{x}_N$ near $\mathfrak{x}$ (if any), the intersection is transversal at each of these points, and the intersection number  of $\mathfrak{H}_1,\ldots,\mathfrak{H}_k$ at $\mathfrak{x}$ equals the sum of the intersection numbers of $\mathfrak{\widetilde H}_1,\ldots,\mathfrak{\widetilde H}_k$  at  $\mathfrak{x}_1,\ldots,\mathfrak{x}_N$.
\end{utver}

From now on, let $\mathbb{K}$ denote either $\C$ or $\T$.

Let $A\subset\Z^n$ be a finite set containing 0. Let $\mathbb{K}^A$ be the set of Laurent polynomials of the form $\sum_{a\in A} c_ax^a,\, c_a\in\mathbb{K}$, and let its subset $(\mathbb{K}\setminus 0)_1^A$ consist of all polynomials, whose constant term $c_0$ equals 1, and whose other coefficients are non-zero.
Let $S\subset(\mathbb{K}\setminus 0)_1^A$ be the $(S_1,\ldots,S_k)$-multisingularity stratum, and let $B_{g_i}\subset(\mathbb{K}\setminus 0)_1^A,\, g_i\in G_i,$ be a family of conditions of incidence parameterized by the set $G_i$. For instance, $g_i\in G_i=(\mathbb{K}\setminus 0)^n$, and $B_ {g_i}$ is the set of all $f\in(\mathbb{K}\setminus 0)_1^A$,
such that $f(g_i)=0$; or $g_i\in G_i=\C^\Delta,\, \Delta\subset\Z^n$, and $B_{g_i}$ is the closure of the set of all $f$ such that $f=0$ is tangent to $g_i=0$ at some point.

We shall be interested in counting hypersurfaces of the form $f=0,\, f\in\CC_1^A$, with prescribed singularities $S_1,\ldots,S_k$, satisfying prescribed conditions of incidence $B_{g_1},\ldots,B_{g_I}$ for generic parameters $g_i\in G_i$. More accurately, we shall be interested in the number of {\it isolated} intersection points of the stratum $S$ and the pencils $B_{g_i}$ for generic $g_i$. Most of the classical problems of enumerative geometry can be restated in this form.
\begin{rem} However, note that the desired number is not in general equal to the global intersection number of $S$ and $B_{g_i}$ in the ring of conditions, because the intersection $S\cap\bigcap_i B_{g_i}$ may have non-isolated components even for generic $g_i$. A classical example is the count of conics that are tangent to five generic lines: the intersection of the corresponding pencils $\cap_{i=1}^5 B_{g_i}$ in the 5-dimensional space of conics contains a two-dimensional component that consists of all two-fold lines. \end{rem}

We now describe the tropical version of counting hypersurfaces with prescribed singularities and conditions of incidence.
Let $\mathfrak{S}$ be the tropical fan of the stratum $S$. It is a tropical fan in the real part of the Lie algebra of the complex torus $\CC_1^A$, this real part is naturally identified with the tropical torus $\TT_1^A$.
\begin{defin} A tropical hypersurface $f=0$ is said to have the $(S_1,\ldots,S_k)$-multisingularity, if $f$ is contained in $\mathfrak{S}$.
 \end{defin}
Assume that each parameter space $G_i$ is a complex torus, denote the real part of its Lie algebra by $\mathfrak{G}_i$, and consider the tropical fan $\mathfrak{B}_i$ of the set $\{(f,g)\, |\, f\in B_g,\, g\in G_i\}$. 
\begin{defin} For every $\mathfrak{g}_i\in \mathfrak{G}_i$, define the {\bf tropical condition of incidence} $\mathfrak{B}_{\mathfrak{g}_i}$ as the tropical intersection of $\mathfrak{B}_i$ and the plane $\TT_1^A\times\{\mathfrak{g}_i\}$. \end{defin}
It is convenient to give tropical conditions of incidence a more geometric interpretation in important special cases. For instance, if $g_i\in G_i=\CC^n$, and $B_{g_i}$ is the set of all $f\in\CC_1^A$, such that $f(g_i)=0$, then $\mathfrak{g}_i\in \mathfrak{G}_i=\TT^n$, and $\mathfrak{B}_{g_i}$ is the set of all $\mathfrak{f}\in\TT_1^A$, such that $\mathfrak{f}(\mathfrak{g}_i)=\mathfrak{0}$. If $g_i\in G_i=\C^\Delta$, $\Delta\in\Z^2$, and $B_{g_i}$ is the closure of the set of all $f\in\CC_1^A$, such that the curves $f=0$ and $g=0$ are tangent, then $\mathfrak{g}_i\in \mathfrak{G}_i=\TT_1^A$, and $\mathfrak{B}_{g_i}$ is the the set of all $\mathfrak{f}\in\TT_1^A$, such that the tropical curves $\mathfrak{f}=0$ and $\mathfrak{g}_i=0$ are tangent in the sense of \cite{MB}.

\begin{defin}
The problem of counting hypersurfaces with prescribed singularities $S_1,\ldots,S_k$, satisfying prescribed conditions of incidence $B_{g_1},\ldots,B_{g_I}$, is said to be {\bf tropicalizable}, if the answer equals the number of isolated points in the intersection $\mathfrak{S}\cap\bigcap_i \mathfrak{B}_{\mathfrak{g}_i}$ for generic $\mathfrak{g}_i\in L_i^*$ (counted with intersection multiplicities, see Definition \ref{deftropint1}).
\end{defin}
In these terms, the correspondence theorem \cite{M} reduces to the computation of the tropical fan of the Severi variety, because the problem of counting plane curves with a prescribed number of self-intersections through prescribed generic points $\mathfrak{g}_1,\ldots,\mathfrak{g}_I$, is tropicalizable (see Proposition \ref{tropinc0} below). Recall that the Severi variety $S$ is the closure of the set of all $f\in\CC^A,\, A\subset\Z^2$, such that the curve $f=0$ is nodal with $d$ nodes.

Similarly, the correspondence theorem \cite{MB} reduces to computation of the tropical fan $\mathfrak{S}$ of the Severi variety, because the problem of counting plane rational curves, containing prescribed generic points and tangent to prescribed generic lines $\mathfrak{g}_1,\ldots,\mathfrak{g}_I$, is tropicalizable (this is difficult, but can be extracted from \cite{MB}). More specifically, the tropical curves counted in the correspondence theorems are given by the equations $\mathfrak{f}=0$, such that $\mathfrak{f}$ is an isolated point of the intersection $\mathfrak{S}\cap\bigcap_i \mathfrak{B}_{\mathfrak{g}_i}$, and the multiplicity, assigned to such a curve by the correspondence theorem, equals the tropical intersection multiplicity of $\mathfrak{S}$ and $\mathfrak{B}_{\mathfrak{g}_i},\, i=1\ldots,I$, at $\mathfrak{f}$.

However, the class of tropicalizable enumerative problems is much wider than counting nodal curves. 
E.g. the problem of counting hypersurfaces with an arbitrary given multisingularity, passing through a given collection of generic points, is tropicalizable:
\begin{utver} \label{tropinc0} Let $S$ be any codimension $p$ subset of $\CC_1^A$ (for example, a multisingularity stratum). Then, for generic points $g_1,\ldots,g_p\in\CC^n$ and $\mathfrak{g}_1,\ldots,\mathfrak{g}_p\in\TT^n$, the intersection number of $S$ and the conditions of incidence $\{f\in\CC_1^A\, |\, f(g_i)=0\}$ coincides with the tropical intersection number of the tropical fan $\mathfrak{S}$ and the conditions of incidence $\{\mathfrak{f}\in\TT_1^A \, |\, \mathfrak{f}(\mathfrak{g}_i)=\mathfrak{0}\},\, i=1,\ldots,p$.
\end{utver}
{\sc Proof.} 
The condition of passing through a generic point has no base points at infinity by Definition \ref{defbpinf}, so the sought equality follows from Proposition \ref{bpinf}(2).  $\quad\square$
\begin{rem}
The intersection number of $S$ and the conditions of incidence $\{f\in\CC_1^A\, |\, f(g_i)=0\}$ equals the number of $f\in S$ such that $f(g_1)=\ldots=f(g_p)=0$, counted with their local intersection multiplicities. Conjecturally, all of these local multiplicites are equal to 1 (i.e. generic conditions of incidence are transversal to $S$), provided that $S$ is a multisingularity stratum and $n>1$. Although this conjecture is well known to be true for $(\mathcal{A}_1,\ldots,\mathcal{A}_1)$--multisingularity strata (the Severi varieties), even in this case it may seem surprising. For instance,  it is not valid for $n=1$, and for any $n$ the $\mathcal{A}_1$--singularity stratum (the $A$-discriminant) is the envelope of the conditions of incidence (i.e. is tangent to the hyperplane $\{f\in\CC_1^A\, |\, f(g)=0\}$ for every $g$).
\end{rem}

\subsection{Tropical fan of the discriminant} \label{Sdiscr}

For a finite set $A$ in the character lattice $L$ of the complex torus $\CC^n$, let us define the {\bf tautological polynomial} $s$ on $\CC^n\times\CC^A$ by the equality $s(x,f)=f(x)$. 
We shall describe the tropical characteristic classes of the critical complete intersection $s=\partial s/\partial x_1=\ldots=\partial s/\partial x_n=0$, where $(x_1,\ldots,x_n)$ are the standard coordinates on $\CC^n$. As a corollary, we shall obtain a new description of the Newton polytope of the $A$-discriminant, provided that it is a hypersurface, and a new characterisation of dual defective toric varieties (projective toric varieties such that their projectively dual variety is not a hypersurface).

The tropical fan of $s=\partial s/\partial x_\cdot=0$ lives in the dual character lattice of $\CC^n\times\CC^A$, i.e. in $L^*\oplus\TT^A$. We shall use the following notation to describe it. For every $\gamma\in L^*\oplus\TT^A$, consider its components $\gamma'\in L^*$ and $\gamma''\in\TT^A,\, \gamma=\gamma'+\gamma''$, as functions on $A$: the function $\gamma'$ is the restriction of the function $\gamma':L\to\Z$ to $A\subset L$, and the value of $\gamma''$ at $a\in A$ is the coefficient of the degree $a$ monomial in the tropical polynomial $\gamma''\in\TT^A$. For every sequence of sets $A_k\subset L$, denote the {affine span} of $\bigcup_{i<k}A_k$ by $A_{<k}$. Assume that $A$ is not contained in an affine hyperplane.
\begin{defin} 1) For every $\gamma\in L^*\oplus\TT^A$, define a finite sequence of subsets $A^\gamma_k\subset A$ inductively by $k\geqslant 0$ as follows: the set $A^\gamma_k$ consists of all points at which the restriction of the function $\gamma''-\gamma':A\to\Z$ to $A\setminus A^\gamma_{<k}$ attains its maximal value, provided that $A^\gamma_{<k}\subsetneq L$, and $A^\gamma_k$ is not defined otherwise.

2) For a finite sequence of non-empty non-intersecting subsets $A_k\subset A$, the cone $C_{A_\bullet}$ is the set of all $\gamma\in L^*\oplus\TT^A$ such that the sequence $A^{\gamma}_k$ equals the sequence $A_k$. The fan that consists of the cones $C_{A_\bullet}$, as $A_\bullet$ runs over all sequences of non-empty non-intersecting subsets of $A$, will be denoted by $C_A$.
\end{defin}
\begin{rem} The codimension of the cone $C_{A_\bullet}$ equals $\sum_i (|A_i|-1)$.
\end{rem}
\begin{rem} It would be useful to prove that the fan $C_A$ is regular (i.e. can be represented as the set of domains of linearity of a convex piecewise-linear function).
\end{rem}

For every sequence of sets $A_k\subset L$, decompose the space $L^\star$ (see Section \ref{Scci}) into a direct sum $\bigoplus_i A^k$ such that $\bigoplus_{i\geqslant k} A^i$ is the orthogonal complement to $A_{<k}$ (that is, $\{l\,|\, l(A_{<k})=0\}$) for every $k$. For every subset $B\subset A$ and a polynomial $f(x)=\sum_{a\in A}c_ax^a$, denote the {\bf restriction of $f$ to $B$} by $f^B(x)=\sum_{a\in B}c_ax^a$, and define the corresponding tautological polynomial $s^B$ on $\CC^n\times\CC^A$   by the equality $s^B(f,x)=f^B(x)$.
Applying Proposition \ref{compcci}(2) to the critical complete intersection $s=\partial s/\partial x_\cdot=0$, we obtain the following.
\begin{utver} \label{comptaut} 1) The closure of $s=\partial s/\partial x_\cdot=0$ in the $C_A$-toric variety intersects its orbits transversally. 

2) The closure of $s=\partial s/\partial x_\cdot=0$ intersects the $C_{A_\bullet}$-orbit by the same set as the closure of the complete intersection given by the equations $\partial_{A^i}s^{A_i}=0$ for all $i$.
\end{utver}

This fact leads to the following description of the characteristic classes of $s=\partial s/\partial x_\cdot=0$.

\begin{defin} 
The {\bf rank} $\rk_k B$ of a subset $B\subset A_k$ is $|B|-\dim($affine span of $B\cup A_{<k})+\dim A_{<k}-1$. A 
sequence $A_k$ is said to be {\bf essential}, if $rk_k B<rk_k A_k$ for every $k$ and every $B\subsetneq A_k$.  
The number of minimal subsets $B\subset A_k$, such that $rk_k B=rk_k A_k-1$, is denoted by $m(A_k)$. The number $m(A_k)$ for the maximal $k$, such that $\rk_k A_k>0$, is denoted by $m(A_\bullet)$. If the number $\bigl|L\bigl/\bigr.\bigcup_i \{a-b\,|\, a \mbox{ and } b\in A_i\}\bigr|$ is finite, then denote it by $i(A_\bullet)$, otherwise set $i(A_\bullet)=0$.
\end{defin}

\begin{exa} We shall work out all the notation for the following examples (in what follows, $A_0$ is always depicted in solid lines and $A_1$ in dashed lines).

\begin{center}
\noindent\includegraphics[width=8.5cm]{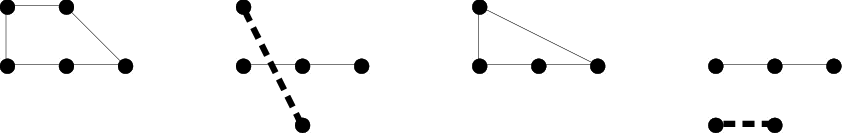}
\end{center}

In the first example, $i(A_\bullet)=1,\, D_0=L^\star$, $\codim C_{A_\bullet}=|A_0|-1=4$, $\rk_0 A_0=2$, and the four minimal sets $B\subset A_0$ with rank $1$ are as follows (so that $m(A_\bullet)=4$):

\begin{center}
\noindent\includegraphics[width=8.5cm]{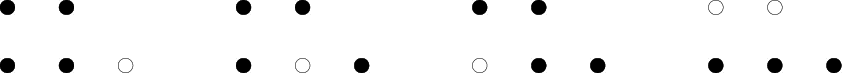}
\end{center}

In the second example, $i(A_\bullet)=2,\, \codim C_{A_\bullet}=|A_0|-1+|A_1|-1=3$, $D_1$ is generated by the function $l(x,y)=x$, and $D_0$ is any plane such that $D_0+D_1=L^\star$.
The third example is 
not essential, and the last one is essential with $i(A_\bullet)=0$.

\end{exa}

\begin{sledst} \label{chartaut} 1) The support set of the characteristic class $\langle s=\partial s/\partial x_\cdot=0\rangle_j$ is contained in the union of the cones $C_{A_\bullet}$ where $A_\bullet$ runs over all essential sequences such that $\sum_i (|A_i|-1)=j$.

2) For $j\leqslant n$, we have $\langle s=\partial s/\partial x_\cdot=0\rangle_j=0$.

3) The multiplicity of the cone $C_{A_\bullet}$ in the tropical fan $\langle s=\partial s/\partial x_\cdot=0\rangle_{n+1}$ equals 1 for every essential $A_\bullet$ such that $\sum_i (|A_i|-1)=n+1$. 

4) The multiplicity of the cone $C_{A_\bullet}$ in the tropical fan $\langle s=\partial s/\partial x_\cdot=0\rangle_{n+2}$ equals $2-m(A_\bullet)$ for every essential $A_\bullet$ such that $\sum_i (|A_i|-1)=n+2$. 
\end{sledst}
{\sc Proof.} 
Choose any sequence $A_\bullet$ and denote $\sum_i (|A_i|-1)$ by $d$. Proposition \ref{comptaut}.2 implies that the intersection of the closure of $s=\partial s/\partial x_\cdot=0$ and the $C_{A_\bullet}$-orbit in the $C_A$-toric variety is a $(d-n-1)$-dimensional variety $V$ in a $d$-dimensional complex torus. The following facts immediately follow from the definitions:

1) If $A_\bullet$ is 
not essential, then $V$ is empty.

2) If $A_\bullet$ is essential, and $d=n+1$, then $V$ consists of one point.

3) If $A_\bullet$ is essential, and $d=n+2$, then $V$ is a projective line minus $m(A_\bullet)$ points.

The statement follows from these facts and Propositions \ref{comptaut} and \ref{schoncc}. $\quad\square$

One can compute multiplicities of all cones for all the characteristic classes of $s=\partial s/\partial x_\cdot=0$ in the same way as we do in the proof of Parts 3 and 4 for the first two non-zero classes. 
The complete answer will be given elsewhere, because we do not need it to study multisingularities in codimension 1 and 2.
\begin{exa} It may well happen that the multiplicity $2-m(A_\bullet)$ is 0 in the circumstances of Part 4, i.e. the cone  $C_{A_\bullet}$ for some essential $A_\bullet$ with $\sum_i (|A_i|-1)=n+2$ does not contribute to $\langle s=\partial s/\partial x_\cdot=0\rangle_{n+2}$. For example, this happens for $A_0=\{(0,0,0),\,(0,0,\pm 1)\}$ and $A_1=\{(0,\pm 1,0),\,(\pm 1,0,0)\}$. In particular, Part 1 gives a non-sharp estimate on the support sets of the characteristic classes. 
\end{exa}

The aim of our computation is the image $\mathfrak{\pi}_*\langle s=\partial s/\partial x_\cdot=0\rangle_j$ under the projection $\pi:\CC^n\oplus\CC^A\to\CC^A$ 
rather than the class $\langle s=\partial s/\partial x_\cdot=0\rangle_j$ itself. 
Let $w_{A_\bullet}$ be the multiplicity 
of the cone $C_{A_\bullet}$ in the characteristic class $\langle s=\partial s/\partial x_\cdot=0\rangle_j,\, \sum_i (|A_i|-1)=j$. Then, by Proposition \ref{projtrop2} and Corollary \ref{chartaut}.1, we have the following.
\begin{utver} \label{imds}
The tropical fan $\mathfrak{\pi}_*\langle s=\partial s/\partial x_\cdot=0\rangle_j$ consists of the cones $\pi(C_{A_\bullet})$ with multiplicities 
$w_{A_\bullet}\cdot i(A_\bullet)$, where $A_\bullet$ runs over all essential sequences such that $\sum_i (|A_i|-1)=j$ and $i(A_\bullet)>0$.
\end{utver}

Together with Corollary \ref{chartaut}.2---4, this completely describes $\pi_*\langle s=\partial s/\partial x_\cdot=0\rangle_j$ for $j\leqslant n+2$. In particular, since the closure of the image of the complete intersection $s=\partial s/\partial x_\cdot=0$ under the projection $\pi:\CC^n\times\CC^A\to\CC^A$ is the $A$-discriminant, 
this gives new information about tropical fans of $A$-discriminants.
Recall that a finite set in $\Z^n$ is called an {\bf iterated circuit}, if it can be decomposed into the union of an essential sequence of non-empty non-intersecting subsets $A_k$, such that $\sum_i (|A_i|-1)=n+1$ and $\dim\sum_i A_i=n$. A more explicit version of this definition is given in the introduction and in \cite{Edcg}, where Part 2 of the following statement was conjectured. 
\begin{sledst} \label{newtdiscr} 1) The dual fan of the Newton polytope of the $A$-discriminant equals $\pi_*\langle s=\partial s/\partial x_\cdot=0\rangle_{n+1}$.

2) The $A$-toric variety is not dual defective if and only if $A$ contains an iterated circuit.
\end{sledst}
{\sc Proof.} Part 1 follows from Remark \ref{rcondfunct}. 
The ``if'' direction of Part 2 is proved in \cite{Edcg}, so it remains to prove the ``only if'' direction: 
if the $A$-toric variety is not dual defective, then the $A$-discriminant is a non-empty hypersurface, then the dual fan of its Newton polytope is not empty, then $\langle s=\partial s/\partial x_\cdot=0\rangle_{n+1}\ne 0$ according to Part 1, then, by Corollary \ref{chartaut}(1), there exists an essential sequence of subsets $A_k\subset A$ such that $\sum_i (|A_i|-1)=n+1$. Its union is an iterated circuit. $\quad\square$

Actually the ``if'' direction of Part 2 also follows from Corollary \ref{chartaut}(3) and Proposition \ref{imds}, because, for every sequence $A_\bullet$, such that $\dim\sum_i A_i=n$, we have $i(A_\bullet)>0$.

\subsection{Tropical fan of the $(\times\times)$ and $(\prec)$ strata: the assumptions} \label{Smain0}

We now 
compute tropical fans of the codimension 2 multisingularity strata. Let $A$ be a finite set in  the character lattice $L$ of the complex torus $\CC^n$.

Let $\{\times\}$, $\{\times\times\}$ and $\{\prec\}$ be the sets of all $f\in\C^A$, such that the closure of the hypersurface $\{f=0\}$ in the smooth toric variety $X_A$ has one ${\mathcal A}_1$ singularity, 
two ${\mathcal A}_1$ singularities, and one ${\mathcal A}_2$ singularity 
respectively. Denote the tropical characteristic classes of these strata by $\langle \times\rangle$, $\langle \times\times\rangle$ and $\langle \prec\rangle$. Note that, unlike in the classical multisingularity theory, we do not switch to the closures of the sets $\{\times\}$, $\{\times\times\}$ and $\{\prec\}$: these sets are smooth (which is not important for us), but not closed! Passing to the closure would not change their tropical fans (that is, their highest characteristic class), but would affect the other characteristic classes.

Theorems \ref{th1}, \ref{th2} and \ref{th3} below provide three independent linear equations on the tropical fans $\langle \times\rangle_2$, $\langle \times, \times\rangle_2$ and $\langle \prec\rangle_2$. We can solve this linear system and compute $\langle \times\times\rangle_2$ and $\langle \prec\rangle_2$, the desired tropical fans  of the $\{\times\times\}$ and $\{\prec\}$ multisingularity strata. In Section \ref{Sexample}, we give an example of the explicit computation of $\langle \times\rangle_2$, $\langle \times\times\rangle_2$ and $\langle \prec\rangle_2$ for a particular set $A$ (i.e. we list all cones and their multiplicities 
in these tropical fans).

We shall obtain these results under certain assumptions on $A$, resembling non-dual defectiveness.

\begin{defin} \label{defvers} For a Laurent polynomial $g\in \C[L]=\C[x_1,x_1^{-1},\ldots,x_n,x_n^{-1}]$, let $I_g\subset\C[L]$ be the ideal generated by $g$ and its partial derivatives $\partial g/\partial x_i$.
The set $A$ is said to be versal in codimension $k$, if the set of all points $f\in\C^A$, such that the natural map $\C^A\to \C[L]/I_f$ is not of the full rank at $f$, has codimension $>k$. 
\end{defin}

Informally, this means that, for each tuple of singularities $S_1,\ldots,S_m$, whose Tyurina numbers sum up to $M\leqslant k$, the $(S_1,\ldots,S_m)$-multisingularity stratum has codimension $M$ and the same adjacencies as in the product of the versal deformations of $S_1,\ldots,S_m$ (c.f. the notions of T-smoothness and deformation completeness of families of curves, see e.g. \cite{gls2}). For codimension 2, we can restate this definition more explicitly:

\begin{defin} \label{propvers}
The set $A$ is versal in codimension 2, if there exists a codimension 3 set $\Sigma\subset\C^A$ such that every $f\in\C^A\setminus\Sigma$ is among the following three types:

1) $f$ is not in the image of the projection $\{s=\partial s/\partial x_\cdot=0\}\to\CC^A$. In this case, the hypersurface $f=0$ has no singularities.

2) $f$ is a regular value of the projection $\{s=\partial s/\partial x_\cdot=0\}\to\CC^A$ with one or two preimages $(z_{(i)},f)$. In this case, $z_{(i)}$ are the only singular points of the hypersurface $f=0$, and both of them are ${\mathcal A}_1$. Moreover, the tautological hypersurface $s=0$ is given by the equation $z_{(i)1}^2+\ldots+z_{(i)n}^2=y_i$ for suitable local coordinates $(y_1,\ldots,y_{|A|})$ near $f\in\C^A$ and $(z_{(i)1},\ldots,z_{(i)n},y_1,\ldots,y_{|A|})$ near $(z_{(i)},f)$.

3) $f$ is a critical value of the projection $\{s=\partial s/\partial x_\cdot=0\}\to\CC^A$ with one preimage $(z,f)$. In this case, $z$ is the only singular point of the hypersurface   $f=0$, and it is ${\mathcal A}_2$. Moreover, the tautological hypersurface $s=0$ is given by the equation $z_1^3+y_1z_1+z_2^2\ldots+z_n^2=y_2$ for suitable local coordinates $(y_1,\ldots,y_{|A|})$ near $f\in\C^A$ and $(z_1,\ldots,z_n,y_1,\ldots,y_{|A|})$ near $(z,f)$.
\end{defin}

In what follows, we always assume $A$ to be versal in codimension 2. For instance, every non-dual defective set in $\Z^2$ and every set containing $4\cdot($simplex of volume 1) is versal. 
\begin{rem} It would be interesting to obtain more accurate sufficient conditions of versality in small codimension and prove its mototonicity similarly to the monotonicity of non-dual defectiveness (if a subset of a finite set $A$ in the lattice $L$ is not dual defective and is not contained in a proper affine sublattice of $L$, then $A$ is not dual defective as well, see \cite{Edcg}). One obvious but weak version of monotonicity is as follows: $A\subset L$ will be called {\it $k$-strongly versal}, if it is versal, and every polynomial $g\in\C[L]$, such that $\dim_\C \C[L]/I_g\leqslant k$, admits $f\in\C^A$, such that $I_g=I_f$. If a $k$-strongly versal $A$ is not contained in a proper affine sublattice of $L$ and is contained in $B$, then $B$ is $(k-1)$-strongly versal.
\end{rem}

In what follows, we also always assume that $A\subset L$ is simple in edges (\cite{vladlen}), i.e. $A$ is the set of lattice points in a lattice polytope, such that every edge is contained in $n-1$ facets, and the $n-1$ external normal covectors to these facets form a part of a basis in $(\Z^n)^*$. Geometrically, this means that the $A$-toric variety $X_A$ has at most isolated singularities.
This assumption is similar to versality of $A$ in codimension 2 with respect to singularities at infinity, see the proposition below. 

For every $B\subset\Z^n$ and every $f=\sum_{a\in A}c_ax^a\in\CC^A$, let $f^B$ be the polynomial $\sum_{a\in B}c_ax^a$, and denote $s^B(x,f)=s(x,f^B)$. If $B$ is a face of the convex hull of $A$, then let $T_B\subset\CC^n$ be the subtorus given by the equations $x^a=x^b$ for all $a$ and $b\in B$, so that the $B$-orbit of the toric variety $X_A$ is naturally identified with $\CC^n/T_B$. Similarly, the intersection of $(B$-orbit$)\times\CC^A$ with the closure of $\{s^B=0\}$ is naturally identified with $\{s^B=0\}/T_B$.
For every facet $\Gamma$ of the convex hull of $A$, let $\Gamma'$ be the maximal face of the convex hull of $A\setminus\Gamma$, parallel to $\Gamma$ and ``looking in the same direction'' (so that $\Gamma'+\Gamma$ is a face of $(A\setminus\Gamma)+A$).  

\begin{utver} \label{lverysm} 
Assume that $A$ is 
simple in edges. 
Then one can choose a codimension 3 set $\Sigma\subset\C^A$ such that every $f\in\C^A\setminus\Sigma$ is among the following two types:

1) $f$ is in the image of the projection $\{s^\Gamma=\partial s^\Gamma/\partial x_\cdot=s^{\Gamma'}=0\}/T_\Gamma\to\CC^A$ for some facet $\Gamma$ of the convex hull of $A$. In this case, the unique preimage of $f$ is a point $(z,f)$, where $z\in(\Gamma$-orbit of $X_A)$ is the unique singular point of the closure of $f=0$, and this singularity is ${\mathcal A}_1$. Moreover, $(\Gamma$-orbit$)\times\CC^A$ and the closure of the hypersurface $s=0$ are given by the equations $z_1=0$ and $(z_1-y_1)^2+z_2^2+\ldots+z_n^2=y_2$ respectively  for suitable local coordinates $(y_1,\ldots,y_{|A|})$ near $f\in\C^A$ and $(z_1,\ldots,z_n,y_1,\ldots,y_{|A|})$ near $(z,f)$.

2) Otherwise, the closure of $f=0$ in the toric variety $X_A$ is smooth at every point outside $\CC^n$.
\end{utver}
The proof is based on the same arguments as Proposition \ref{compcci} (note that the lattice distance between the parallel hyperplanes, containing $\Gamma$ and $\Gamma'$, equals 1, because $A$ is simple in edges).

\subsection{Tropical fan of the $(\times\times)$ and $(\prec)$ strata: the answer} \label{Smain}

Recall that we assume that $A\subset\Z^n$ is simple in edges and versal in codimension 2 (see Section \ref{Smain0}).

\vspace{1ex}

{\sc The first equation.} 

\begin{theor} \label{th1} If $A$ is versal in codimension 2 and simple in edges, then 

$$\langle \times\rangle_2+2\langle \times\times\rangle_2+2\langle \prec\rangle_2=\sum_{\Gamma} \pi_*\langle \{s^{\Gamma}=0\}/T_{\Gamma}\rangle_{\dim\Gamma+2},$$

where $\Gamma$ runs over all positive-dimensional faces of the convex hull of $A$.
\end{theor}

Here $\langle \{s^{\Gamma}=0\}/T_{\Gamma}\rangle_{j}$ for all $j$ can be expressed in terms of the Newton polytope of $s^{\Gamma}$, see Example \ref{trophyp}. So the only unknown terms are on the left hand side.

\vspace{1ex} 

{\sc Proof.} The sought equality is the statement of Proposition \ref{eqstar} for $i=2$ and the projection $\pi:X_A\times\CC^A\to\CC^A$, restricted to the closure of the tautological hypersurface $\{s=0\}\subset X_A\times\CC^A$, $s(x,f)=f(x)$.  More complicated multisingularities do not appear in the left hand side by the versality of $A$ in codimension 2 and by Proposition \ref{lverysm}. $\quad\Box$

\vspace{1ex} 

{\sc The second equation.} 

\begin{theor} \label{th2} If $A$ is versal in codimension 2 and simple in edges, then 

$$\langle \times\rangle_2+2\langle \times\times\rangle_2+\langle \prec\rangle_2=\pi_*\langle s=\frac{\partial s}{\partial x_{\cdot}}=0\rangle_{n+2}+\sum_{\Gamma} \pi_*\langle \{s^{\Gamma}=s^{\Gamma'}=\frac{\partial s^{\Gamma}}{\partial x_{\cdot}}=0\}/T_{\Gamma}\rangle_{n+1},$$

where $\Gamma$ runs over all facets of the convex hull of $A$.
\end{theor}

Here $\langle \{s^{\Gamma}=s^{\Gamma'}=\frac{\partial s^{\Gamma}}{\partial x_{\cdot}}=0\}\rangle_{n+1}=\langle \{s^{\Gamma}=\frac{\partial s^{\Gamma}}{\partial x_{\cdot}}=0\}\rangle_{n}\cdot\langle s^{\Gamma'}=0\rangle_1$ and $\langle s=\frac{\partial s}{\partial x_{\cdot}}=0\rangle_{n+2}$ are computed by Corollary \ref{chartaut}. So the only unknown terms are on the left hand side.

\vspace{1ex}

{\sc Proof.} The sought equality is the statement of Proposition \ref{eqstar} for $i=2$ and the projection $\pi:X_A\times\CC^A\to\CC^A$, restricted to the closure of the first Boardman stratum $\{s=\frac{\partial s}{\partial x_{\cdot}}=0\}$. More complicated multisingularities do not appear in the left hand side by the versality of $A$ in codimension 2 and by Proposition \ref{lverysm}. The explicit expression for the intersection with the orbit in the right hand side is given by Proposition \ref{comptaut}(2). $\quad\Box$

\begin{rem} The equalities from Theorems \ref{th1} and \ref{th2} are already enough to compute the tropical fan of the singularity stratum $\{\prec\}$. In the same way as we deduced Theorems \ref{th1} and \ref{th2}, we can also compute the left hand side of the equality in Proposition \ref{eqstar} for all higher Boardman strata of the projection $\{s=0\}\to\C^A$. However, for the $(1,1)$-Boardman stratum (i.e. the set of critical points of the restriction of $\pi$ to its set of critical points) we shall obtain an equation on $\langle \prec\rangle_2$ that is linearly dependent with the equalities of Theorems \ref{th1} and \ref{th2}, and for higher Boardman strata (see e.g. \cite{avg}) the equations will be trivial in codimension 2. So Proposition \ref{eqstar} alone is not enough to find the tropical fan even for the simplest multisingularity stratum $\{\times\times\}$. \end{rem}

\vspace{1ex} 

{\sc The third equation.} To formulate it, we define a certain continuous piecewise-linear function $\mathbf{l}_A$ on the support set of the tropical fan $\langle s=\frac{\partial s}{\partial x_{\cdot}}=0\rangle_{n+1}$. Its definition relies upon the notions and notation from Section \ref{Sdiscr}.

Consider $\gamma\in \TT^n\oplus\TT^A$  such that $i(A^\gamma_\bullet)>0$. Recall that we consider its components $\gamma'\in \TT^n$ and $\gamma''\in\TT^A,\, \gamma=\gamma'+\gamma''$, as functions on $A$: the function $\gamma'$ is the restriction of the function $\gamma':\R^n\to\Z$ to $A\subset\Z^n$, and the value of $\gamma''$ at $a\in A$ is the coefficient of the degree $a$ monomial in the tropical polynomial $\gamma''\in\TT^A$. For every $k$, denote the value of the function $\tilde\gamma=\gamma''-\gamma':A\to\Z$ at $A^\gamma_k$ by $\tilde\gamma_k$, and define
$A^\gamma_k(r)$ as the set of all $a\in A$ such that $\tilde\gamma(a)\in[\tilde\gamma_k,\, \tilde\gamma_k+r]$.
Define the function $$i_\gamma(r)=\bigl|\Z^n\bigl/\bigr.\bigcup_i\{a-b\,|\, a \mbox{ and } b\in A^\gamma_k(r)\bigr\}\bigr|$$ on $\R$ and the function $$\mathbf{l}_A(\gamma)=\int_0^{+\infty} (i_\gamma(r)-1)dr$$ on the set $\{\gamma\, |\, i(A^\gamma_\bullet)>0\}$. 

The function $\mathbf{l}_A:\{\gamma\, |\, i(A^\gamma_\bullet)>0\}\to\R$ is well-defined, because $i(A^\gamma_\bullet)>0$ implies
$i_\gamma(r)\leqslant i_\gamma(0)<\infty$, and the assumption that $A$ affinely generates $\Z^n$ implies $i_\gamma(r)=1$ for $r$ large enough. Moreover, the function $\mathbf{l}_A:\{\gamma\, |\, i(A^\gamma_\bullet)>0\}\to\R$ is continuous, because, for $\gamma_1$ in a small neighborhood of $\gamma$, the functions $i_\gamma$ and $i_{\gamma_1}$ coincide outside a small neighborhood of the discontinuities of $i_\gamma$.

Let $\mathbf{\tilde l}_A$ be any continuous extension of the function $\mathbf{l}_A$ to the support set of the fan $\langle s=\frac{\partial s}{\partial x_{\cdot}}=0\rangle_{n+1}$, then the product $\mathbf{\tilde l}_A\cdot \langle s=\frac{\partial s}{\partial x_{\cdot}}=0\rangle_{n+1}$ is an element of the component $\mathcal{K}_{n+1}^1\bigl(\Z^n\oplus\Z^A\bigr)$ of the differential ring of tropical fans (see Section \ref{Stit}). Since its image $\pi_*(\mathbf{\tilde l}_A\cdot \langle s=\frac{\partial s}{\partial x_{\cdot}}=0\rangle_{n+1})$ under the projection $\pi:\CC^n\times\CC^A\to\CC^A$
does not depend on the choice of an extension $\mathbf{\tilde l}_A$ by Proposition \ref{imds}, we shall denote it by $$\pi_*(\mathbf{l}_A\cdot \langle s=\frac{\partial s}{\partial x_{\cdot}}=0\rangle_{n+1})\in \mathcal{K}_{1}^1\bigl(\Z^A\bigr).$$

\begin{defin} \label{defsecf} 
The corner locus $$\delta\pi_*(\mathbf{l}_A\cdot \langle s=\frac{\partial s}{\partial x_{\cdot}}=0\rangle_{n+1}) \in \mathcal{K}_{2}^0\bigl(\Z^A\bigr)$$ is a codimension 2 tropical fan in $\TT^A$ and will be called the {\bf ternary fan} of $A$ (by analogy with the codimension 1 {\bf secondary fan} introduced in \cite{GKZ94} and containing $\pi_*(\langle s=\frac{\partial s}{\partial x_{\cdot}}=0\rangle_{n+1}) \in \mathcal{K}_{1}^0\bigl(\Z^A\bigr)$ as a summand).
\end{defin}

\begin{rem} Note that the corner locus of the product $\mathbf{f}F\in\mathcal{K}_{1}^1$ of a continuous piecewise-linear function $\mathbf{f}\in\mathcal{K}_0^1$ and a tropical fan $F\in\mathcal{K}_1^0$ is known as the Cartier divisor of $\mathbf{f}$ on $F$ (see e. g. \cite{M06}), and the geometric meaning of $\delta\pi_*(\mathbf{l}\langle s=\frac{\partial s}{\partial x_{\cdot}}=0\rangle_{n+1})$ is similar. However, the element $\pi_*(\mathbf{l}\langle s=\frac{\partial s}{\partial x_{\cdot}}=0\rangle_{n+1})\in\mathcal{K}_{1}^1\bigl(\Z^A\bigr)$ that we have to differentiate cannot be decomposed into the product of a continuous piecewise linear function from $\mathcal{K}_{0}^1\bigl(\Z^A\bigr)$ and a tropical fan with constant weights from $\mathcal{K}_{1}^0\bigl(\Z^A\bigr)$.
\end{rem}

\begin{theor} \label{th3} If $A$ is versal in codimension 2 and simple in edges, then
$$ \langle \times\rangle_2 
-\langle \prec\rangle_2=\delta\pi_*(\mathbf{l}\langle s=\frac{\partial s}{\partial x_{\cdot}}=0\rangle_{n+1})-(\pi_*\langle s=\frac{\partial s}{\partial x_{\cdot}}=0\rangle_{n+1})^2.$$
\end{theor}
See Section \ref{Sproofs123} for a proof. The only unknown classes in this theorem are on the left hand side, and, together with Theorems \ref{th1} and \ref{th2}, we obtain three independent equations on $\langle \times\rangle_2$, $\langle \times, \times\rangle_2$ and $\langle \prec\rangle_2$.

\begin{rem} The geometric meaning of the piecewise linear function $\pi_*(\mathbf{l}_A\cdot \langle s=\frac{\partial s}{\partial x_{\cdot}}=0\rangle_{n+1})$ is clear: it is the tropicalization of the function on the $A$-discriminant $D=\pi_*\{s=\frac{\partial s}{\partial x_{\cdot}}=0\}$, assigning the determinant of the Jacobian matrix $\left(\frac{\partial^2f}{\partial x_i\partial x_j}(z)\right)_{i,j=1,\ldots,n}$ to every $f\in D$ with a unique point $z\in\CC^n$ such that $f(z)=0$ and ${\rm d}f(z)=0$. However, I do not know how to prove Theorem \ref{th3} using this interpretation.
\end{rem}

\begin{rem} \label{remth3}
The cones of the fan $(\pi_*\langle s=\frac{\partial s}{\partial x_{\cdot}}=0\rangle_{n+1})^2$ correspond to the two-dimensional faces of the polytope, dual to the fan $\pi_*\langle s=\frac{\partial s}{\partial x_{\cdot}}=0\rangle_{n+1}$. Multiplicities of the cones equal the areas of the corresponding faces. Every such face is either a parallelogram (whose area is easy to compute) or the secondary polygon of a certain subset of $A$ (whose area can be conveniently computed by the formula from Example \ref{secpolygon} below).
\end{rem}

We now conjecture how the tropical fan of the $(\underbrace{\times,\ldots,\times}_m)$-stratum might look like for arbitrary $m$. 
\begin{defin} The tropical fan 
$\delta^k\pi_*(\mathbf{l}^k\langle s=\frac{\partial s}{\partial x_{\cdot}}=0\rangle_{n+1})$
is called the $k$-ary fan of $A$ and is denoted by $A^k$.
\end{defin}
For every linear $\gamma:\Z^n\to\Z$, define $A^\gamma=A_1^\gamma$ as the set of all points 
of $A$ where the restriction $\gamma|_A$ attains its maximum, and then inductively $A^\gamma_{k+1}=(A\setminus\bigcup_{i=1}^k A^\gamma_k)^\gamma$.
\begin{conj} There exists a universal function $P_m$ of a collection of disjoint sets $B_1,\ldots,B_m\subset\Z^n$, taking values in $\mathcal{K}_{m}^0\bigl(\Z^{\bigcup_i B_i}\bigr)$, such that

1) For given $\dim B_1,\ldots,\dim B_m$, the function $P_m(B_1,\ldots,B_m)$ is a polynomial of the $k$-ary fans $(B_i)^k$ and the characteristic classes $\pi_*\langle s_i=\frac{\partial s_i}{\partial x_{\cdot}}=0\rangle_j$ for the tautological polynomials $s_i=\sum_{b\in B_i} c_bx^b$ on $\CC^n\times\CC^{B_i}$;

2) the tropical fan of the $(\underbrace{\times,\ldots,\times}_m)$-stratum in $\CC^A$ is the sum of $P_m(A_1^\gamma,\ldots,A_m^\gamma)$ over all primitive $\gamma$ (including $\gamma=0$).
\end{conj}

\subsection{Reduction to the universal case} \label{Sred}

We have described the tropical fans of the universal multisingularity strata $\{\times\times\}$ and $\{\prec\}$ in $\CC^A$. We now outline how to apply this result to more down-to-earth settings, such as the multisingularity strata $\{\times\times\}$ and $\{\prec\}$ of a generic polynomial map $\CC^m\to\C^n$. The details will be given in a separate paper.

{\sc Non-morse polynomials.} Let $A\not\ni 0$ be a finite set, generating $\Z^n$. The set of polynomials $f\in\CC^A$ such that $f:\CC^n\to\C$ has a degenerate critical point  is called the {\it caustics};  the set of polynomials $f\in\CC^A$ such that $f:\CC^n\to\C$ takes the same value at two critical points is called the {\it Maxwell stratum}, see \cite{landozvon} for the study of the case $n=1$.

Denote $A\cup\{0\}$ by $A'$ and the natural projection $\CC^{A'}\to\CC^{A}$ by $p$, then the tropical fan of the caustics equals $p_*(\langle \prec\rangle_2)$, and the tropical fan of the Maxwell stratum equals $p_*(\langle \times\times\rangle_2)$, where $\{\times\times\}$ and $\{\prec\}$ are the codimension 2 multisingularity strata in $\CC^{A'}$.

{\sc $\{\times\times\}$ and $\{\prec\}$ for complete intersections.} 
For a collection of finite sets $A_0,\ldots,A_k\subset\Z^m$, let $\{\times\}$ be the set of all $(f_1,\ldots,f_k)\in\C^{A_0}\oplus\ldots\oplus\C^{A_k}$ such that $f_1=\ldots=f_k=0$ is a complete intersection with a unique singular point, and its type is $\mathcal{A}_1$. Similarly, define $\{\prec\}$ and $\{\times\times\}\subset\C^{A_0}\oplus\ldots\oplus\C^{A_k}$ as the sets of all complete intersections with one $\mathcal{A}_2$ singularity and two $\mathcal{A}_1$ singularities respectively. Define the Cayley configuration $A_0*\ldots*A_k\subset\Z^n\times\Z^k$ as the union of $A_0\times\{0\}$ and $A_i\times\{i$-th vector of the standard basis$\}$ over $i=1,\ldots,k$, and consider the natural isomorphism 
$p:\CC^{A_0*\ldots*A_k}\to\CC^{A_0}\times\ldots\times\CC^{A_k}$.
\begin{utver} \label{propci} Assume that none of the $A_i$ is contained in an affine hyperplane, and $A_0*\ldots*A_k$ is 2-versal. 
Then we have $$\langle \times\rangle_1=p_*\langle \times\rangle_1,\; \langle \prec\rangle_2=p_*\langle \prec\rangle_2 \mbox{ and } \langle \times\times\rangle_2=p_*\langle \times\times\rangle_2.$$
\end{utver}
{\it Proof.} Denoting the closure of the highest dimension component of a constructible set $S$ by $S_0$, the statement of the proposition follows from the equalities $$\{\times\}_0=p(\{\times\}_0),\, \{\prec\}_0=p(\{\prec\}_0),\, \{\times\times\}_0=p(\{\times\times\}_0).$$
For the discriminant $\{\times\}_0$, this equality is proved in \cite[Theorem 2.31]{Edcg} (see also \cite{Smix}), and for the codimension 2 strata the proof follows the same lines. $\quad\Box$

Note that however the inclusion $p(\{\times\})\subset\{\times\}$ is proper, and thus $\langle \times\rangle_2\ne p_*\langle \times\rangle_2$ in general.

{\sc $\{\times\times\}$ and $\{\prec\}$ for projections of complete intersections.} Let $B_0,\ldots,B_k$ be finite sets in $\Z^n\oplus\Z^m$, $g_i\in\CC^{B_i}$ be generic polynomials, and $\pi$ be the restriction of the projection $\CC^n\times\CC^m\to\CC^m$ to the complete intersection $g_0=\ldots=g_k=0$. We shall compute the tropical fans of the multisingularity strata $\{\times\times\}_\pi$ and $\{\prec\}_\pi\subset\CC^m$ of $\pi$.

For this, denote the image of $B_i$ under the projection $\Z^n\oplus\Z^m\to\Z^n$ by $A_i$, and let $x$ and $y$ be the standard coordinates on $\CC^n$ and $\CC^m$ respectively. 
In the product $\CC^{A_0}\times\ldots\times\CC^{A_k}\times\CC^m$, consider the sets $G_i=\{(f_0,\ldots,f_k,y)\,|\, \forall x\; f_i(x)=g_i(x,y)\}$ and the projections $p$ to $\CC^m$ and $q$ to $\CC^{A_0}\times\ldots\times\CC^{A_k}$. The tropical fans of $\{\times\times\}_\pi$ and $\{\prec\}_\pi$ equal 
 $$p_*( \langle G_0\rangle_{|A_0|}\cdot\ldots\cdot \langle G_k\rangle_{|A_k|}\cdot q^*\langle \times\times\rangle_2) \mbox{ and } p_*(\langle G_0\rangle_{|A_0|}\cdot\ldots\cdot \langle G_k\rangle_{|A_k|}\cdot q^*\langle \prec\rangle_2). \eqno{(*)}$$ In this formula, the tropical fan $\langle G_i\rangle_{|A_i|}$ of the set $G_i$ is known for generic $g_i\in\CC^{B_i}$: the set $G_i$ is a non-degenerate complete intersection, so its tropical fan is the product of the dual fans of the Newton polytopes of its equtions.  Computation of the tropical fans $\langle \times\times\rangle_2$ and $\langle \prec\rangle_2$ on the right hand side is reduced to Theorems \ref{th1}, \ref{th2} and \ref{th3} by Proposition \ref{propci}.
\begin{rem} This result covers the setting of Example \ref{exa1}, but reducing the general answer $(*)$ to the elementary one given in Example \ref{exa1} is a non-trivial combinatorial problem, which will be treated in a subsequent paper. On the other hand, the answer $(*)$ remains valid for generic hypersurfaces $M\subset\CC^3$ with arbitrary Newton polytopes.
\end{rem}

{\sc $\{\times\times\}$ and $\{\prec\}$ for maps of complete intersections.} Let $A_0,\ldots,A_p$, $D_1,\ldots,D_m$ be finite sets in $\Z^n$, $f_i\in\CC^{A_i}$ and $h_i\in\CC^{D_i}$ be generic polynomials, and $H$ be the restriction of the map $(h_1,\ldots,h_m):\CC^n\to\C^m$ to the complete intersection $f_0=\ldots=f_k=0$. We shall compute the tropical fans of the multisingularity strata $\{\times\times\}_H$ and $\{\prec\}_H\subset\CC^m$ of $H$.
For this, note that $\{\times\times\}_H$ and $\{\prec\}_H$ equal $\{\times\times\}_\pi$ and $\{\prec\}_\pi$, where $\pi$ is the restriction of the projection $\CC^n\times\CC^m\to\CC^m$ to the complete intersection $f_0=\ldots=f_p=h_1-y_1=\ldots=h_m-y_m=0$, and $(y_1,\ldots,y_m)$ are the standard coordinates on $\CC^m$. The tropical fans of $\{\times\times\}_\pi$ and $\{\prec\}_\pi$ can now be computed as in $(*)$. 

\begin{rem} If we try to apply Proposition \ref{propci} to computation of the tropical fans $\langle \times\times\rangle_2$ and $\langle \prec\rangle_2$ on the right hand side of $(*)$ in this setting, we observe that the assumption of Proposition \ref{propci} is not satisfied. Nevertheless one can easily verify that the conclusion of Proposition \ref{propci} remains valid in this particular setting.\end{rem}

\subsection{Volume of a fiber body} \label{Svfp}

Choose a codimension $k$ subspace $L\subset\R^n$. There exists a unique additive symmetric function $M_L$ of $k+1$ convex bodies in $\R^n$, taking values in convex bodies in $L$, such that $M_L(\underbrace{A,\ldots,A}_{k+1})$ is the fiber body of $A$ for every convex body $A\subset\R^n$. See Example \ref{exafb} or \cite{bs} for the definition of the fiber body, and \cite{mcmmix} or \cite{mmj08} for a proof of the fact.
\begin{defin}
The convex body $M_L(A_0,\ldots,A_k)$ is called the {\bf mixed fiber body} of $A_0,\ldots,A_k$. \end{defin}
See \cite{EK} for a relation of mixed fiber polytopes to algebraic geometry and the following characterization.
\begin{utver} \label{ekh1} The mixed fiber body of $A_0,\ldots,A_k$ is the unique convex body $X$ such that the Euclidean mixed volume of $A_0,\ldots,A_k,B_1,\ldots,B_{n-k-1}$ in $\R^n$ equals the Euclidean mixed volume of 
$X,B_1,\ldots,B_{n-k-1}$ in $L$.
\end{utver}
\begin{utver} \label{vfp} Let $L$ be a codimension $k$ subspace in $\R^n$, and $A_i^j\subset\R^n,\, i=0,\ldots,k,\, j=1,\ldots,n-k,$ be convex bodies. Denote the space $\{(v_1,\ldots,v_{n-k}\,|\, v_1+\ldots+v_{n-k}=0\}\subset L^{n-k}$ by $K$ and the image of $$\underbrace{\{0\}\times\ldots\times\{0\}}_{i-1}\times A^i_j\times\{0\}\times\ldots\times\{0\}\subset(\R^n)^{n-k}$$ under the projection $(\R^n)^{n-k}\to(\R^n)^{n-k}/K$ by $\widetilde A^i_j$. Then the mixed volume of the mixed fiber bodies $M_L(A_0^i,\ldots,A_k^i),\, i=1,\ldots,n-k$ equals the mixed volume of $\widetilde A_i^j,\, i=0,\ldots,k,\, j=1,\ldots,n-k$.
\end{utver}
The proof is by $n-k-1$ applications of the equality from Proposition \ref{ekh1}.

An important special case of a fiber body is a secondary polytope (\cite{GKZ94}), and Proposition \ref{vfp} gives a formula for its volume (see e.g. \cite{orevk} for one motivation to study volumes of secondary polytopes). 
We shall need the following special case.
\begin{exa} \label{secpolygon}
Let $A\subset\R^n$ be a lattice polytope. For every its vertex $v$, let $c_A^v$ be the volume of the difference of $A$ and the convex hull of $A\cap\Z^n\setminus\{v\}$. If the secondary polytope of $A$ is 2-dimensional, then its lattice area equals $(n+1)!\Vol A-n!\sum_v c_A^v$.
\end{exa}
Another important special case of a mixed fiber body is the Newton polytope of a sparse resultant (\cite{EK}), and Proposition \ref{vfp} gives a positive formula for its volume. This may be already of interest in the case of the Newton polytope of the determinant, that is, the Birkhoff polytope (see e.g. \cite{bir1} for the study of its volume).

\subsection{Proof of the Third equation} \label{Sproofs123}

Let $i=(i_1,i_2,\ldots)$ be a sequence of integer numbers stabilizing at 1, such that $i_r|i_{r-1}$ for every $r$.
\begin{defin} An {\bf$i$-forking paths singularity} is a plane singularity with $i_1$ distinct regular branches $\varphi_{q_1,q_2,\ldots}:(\C,0)\to(\C^2,0),\, q_r=1,\ldots,i_r/i_{r+1}$, such that the intersection number of $\varphi_{p_1,p_2,\ldots}$ and $\varphi_{q_1,q_2,\ldots}$ at 0 equals the minimal number $r$ such that $p_r\ne q_r$. 
\end{defin}

\begin{lemma} \label{milnfork} The Milnor number of an $i$-forking paths singularity equals $\mu(i)=i_1\sum\limits_{r}(i_r-1)-(i_1-1)$.
\end{lemma}
{\sc Proof.} Perturb the branches of the singularity independently, then the union $U$ of the perturbations has exactly $N=\sum_r \frac{i_1}{i_r}\frac{i_r(i_r-1)}{2}$ nodes. The Euler characteristic of the normalization of $U$ equals $i_1$ and differs by $2N$ from the desired Euler characteristic of the Milnor fiber. $\quad\Box$

{\sc Proof of Theorem \ref{th3}.} We can now prove the following equality, equivalent to the statement of Theorem \ref{th3}. Note that the tropical fan $\pi_*\langle s=\frac{\partial s}{\partial x_\cdot}=0\rangle_{n+1}$ is the dual fan of the Newton polytope $\mathcal{N}_A$ of the Gelfand-Kapranov-Zelevinsky $A$-discriminant $D_A$.

\begin{utver} 
For a generic lattice polytope $B\subset\R^A$ (that is, every polytope outside finitely many hypersurfaces in the cone of polytopes), we have
$$ - [\mathcal{N}_A]^2 [B]^{|A|-2} - (|A|-2) [\mathcal{N}_A] [B]^{|A|-1} 
=$$ $$= \langle \times\rangle_2 [B]^{|A|-2} - (|A|-2) \langle\times\rangle_1 [B]^{|A|-1}
-\langle\prec\rangle_2 [B]^{|A|-2} - \delta\pi_*(\mathbf{l}\langle s=\frac{\partial s}{\partial x_{\cdot}}=0\rangle_{n+1}) [B]^{|A|-2}.$$
\end{utver}
This is equivalent to Theorem \ref{th3}, because the second terms on both sides mutually eliminate: $[\mathcal{N}_A]=\langle \times\rangle_1$, and then the multiplier $[B]^{|A|-2}$ can be cancelled by the nondegeneracy of the intersection number pairing in the ring of tropical fans $\mathcal{K}^0(\Z^A)$.

{\sc Proof.} For $M\in\Z$, consider the map $\mathcal{M}:\C^A\to\C^A,\, \mathcal{M}(\sum_a c_ax^a)=\sum_a c_a^{M!}x^a$, rising the standard coordinates to a large power. Choose generic polynomials $h_1,\ldots,$ $h_{|A|-2}\in\C^{\Z^A\cap B}$, generic $g\in\C^{\Z^A\cap M\mathcal{N}_A}$ close to $D_A\circ\mathcal{M}\in\C^{\Z^A\cap M\mathcal{N}_A}$, a smooth toric compactification $\mathcal{X}$ of $\C^A$, compatible with the polytopes $B$ and $\mathcal{N}_A$, and denote its big torus by $T$. Denote the closure of $h_1=\ldots=h_{|A|-2}=0$,  $h_1=\ldots=h_{|A|-2}=g=0$ and $h_1=\ldots=h_{|A|-2}=D_A\circ\mathcal{M}=0$ in $\mathcal{X}$ by $V$, $C$ and $\widetilde{C}$. The smooth curve $\widetilde{C}$ is a smoothening of the curve $C$ in the smooth surface $V$, so we have

$$e(T\cap\widetilde{C}) = e(\mbox{smooth part of }C)+\sum_{p\in{\rm sing}\, C} e\bigl(T\cap(\mbox{Milnor fiber of }C\mbox{ at }p)\bigr)$$
where all points of $C$ outside the big torus $T$ are considered as singular.

The curve $C$ has the following singularities:

1) $\mathcal{A}_2$-singularities appear at the points of the intersection of $V$ and the singularity stratum $\mathcal{M}^{-1}\{\prec\}$ (because the discriminant of the versal deformation of an $\mathcal{A}_2$-singularity is itself a $\mathcal{A}_2$-singularity). 

2) $\mathcal{A}_1$-singularities appear at the points of the intersection of $V$ and the multisingularity stratum $\mathcal{M}^{-1}\{\times\times\}$. 

3) For large $M$, all singularities of $C$ outside the big torus of $\mathcal{X}$ are forking paths singularities (this is not always true for $M=1$).

Thus, the previous formula can be rewritten more precisely as follows:

$$e(T\cap\widetilde{C}) = e(\mbox{smooth part of }C)-(\mbox{the number of }\mathcal{A}_2\mbox{-singularities of }C)+$$ $$+\sum_{p\in T\cap{\rm sing}\, C} e\bigl(T\cap(\mbox{Milnor fiber of }C\mbox{ at }p)\bigr).\eqno{(1)}$$

As we evaluate the terms of this equality, they turn out to be equal to the respective terms of the desired equality in the statement.

By the Khovanskii formula \cite{khov0}, the left hand side of $(1)$ equals $$- M!^2\cdot[\mathcal{N}_A]^2 [B]^{|A|-2} - M!\cdot (|A|-2) [\mathcal{N}_A] [B]^{|A|-1}.$$

For the first term on the right hand side of $(1)$, we have $$e(\mbox{smooth part of }C)=\langle \mbox{smooth part of }C\rangle_{|A|}=(\langle h_1=0\rangle\cdot\ldots\cdot\langle h_{|A|-2}=0\rangle\cdot\langle \times\rangle)_{|A|}=$$ $$=M!^2\cdot\langle \times\rangle_2 [B]^{|A|-2} - M!\cdot (|A|-2) \langle \times\rangle_1 [B]^{|A|-1}.$$

For the second term on the right hand side of $(1)$, the $\mathcal{A}_2$-singularities appear at the points of the intersection of $V$ and the singularity stratum $\mathcal{M}^{-1}\{\prec\}$, so there are $M!^2\cdot\langle \prec\rangle_2 [B]^{|A|-2}$ of them. 

For the last term of $(1)$, let us describe the forking paths singularities of $C$ in more detail.
For a generic lattice polytope $B\subset\Z^n$, the tropical fans $[B]^{|A|-2}$ and $[\mathcal{N}_A]$ intersect transversally (outside the origin), and the 1-dimensional fan $[B]^{|A|-2}[\mathcal{N}_A]$ consists of finitely many rays with multiplicities 
$m_\alpha\in\Z$ and primitive generators $\lambda_\alpha\in(\Z^A)^*$. 
Considering $\lambda_\alpha\in(\Z^A)^*\cong\Z^A$ as a function $A\to\Z$, we can find a linear (but not necessarily integer!) function $\gamma'_\alpha:\Z^n\to\Q$ such that the function $\lambda_\alpha-\gamma'_\alpha:A\to\Z$ attains its maximum on a certain set $A_\alpha\subset A$ that is not the set of vertices of a simplex. Moreover, for a generic $B$, such $A_\alpha$ is unique and is a cirquit. Define $\gamma_\alpha=(\lambda_\alpha,\gamma'_\alpha)\in(\Q^A\oplus\Q^n)^*$. For $r\in\Z$, define $i^\alpha_r$ as $i(\frac{r-1}{M!})$, where the function $i(\cdot)$ is introduced as in Definition \ref{defsecf} for $\gamma=\gamma_\alpha$.

In this notation and under these genericity assumptions for the lattice polytope $B\subset\R^A$, the codimension 1 orbit of $\mathcal{X}$, corresponding to the covector $\lambda_\alpha$, contains $m_\alpha M!/i^\alpha_1$ singularities of the curve $C$, and each of them is an $(i^\alpha_1,i^\alpha_2,\ldots)$-forking paths singularity, all of whose branches are transversal to the orbit. Other orbits of $\mathcal{X}$ do not contain singularities of $C$. Counting the Euler characteristic of the Milnor fibers of the forking paths singularities by Lemma \ref{milnfork}, we conclude that the last term  in $(1)$ equals $-M!^2\cdot \delta\pi_*(\mathbf{l}\langle s=\frac{\partial s}{\partial x_{\cdot}}=0\rangle_{n+1}) [B]^{|A|-2}$. $\quad\Box$

\newpage \subsection{Example} \label{Sexample}

We compute the tropical fans of the $\{\times\times\}$ and $\{\prec\}$-strata in $\C^{A}$, with $A$ as follows:

\begin{center}
\noindent\includegraphics[width=1.7cm]{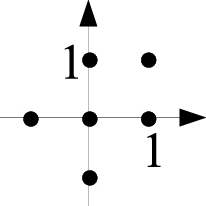}
\end{center}

We shall describe cones in $\Q^A$ that we encounter in the process of computation by drawing coherent subdivisions of the convex hull of $A$: a coherent subdivision will stand for the cone that consists of all $\gamma\in\Q^A$, such that the upper faces $F$ of the convex hull of the graph of $\gamma:A\to\Q$ form the given subdivision. A point $a$ is shown in black or white depending on whether $(\gamma(a), a)$ is in $F$ or not, as in $\cite{GKZ94}$.

The right-hand side in Theorem \ref{th1} consists of the cones

\begin{center}
\noindent\includegraphics[height=1.2cm]{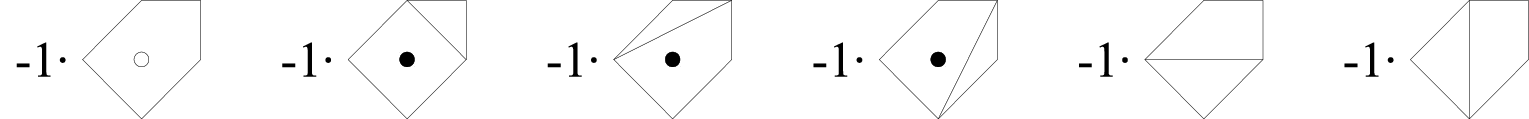}
\end{center}

The right-hand side in Theorem \ref{th2} consists of the cones

\begin{center}
\noindent\includegraphics[height=1.2cm]{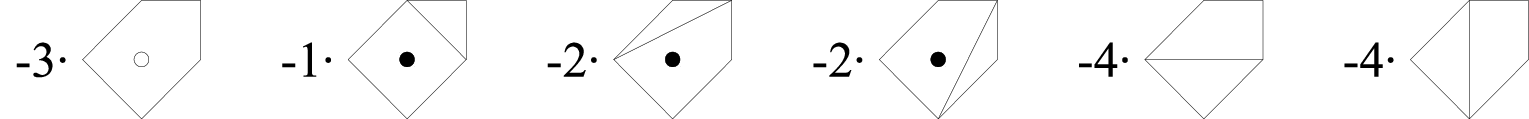}
\end{center}

The fan $d\pi_*(\mathbf{l}\langle s=\frac{\partial s}{\partial x_{\cdot}}=0\rangle_{n+1})$ consists of the cones

\begin{center}
\noindent\includegraphics[height=1.2cm]{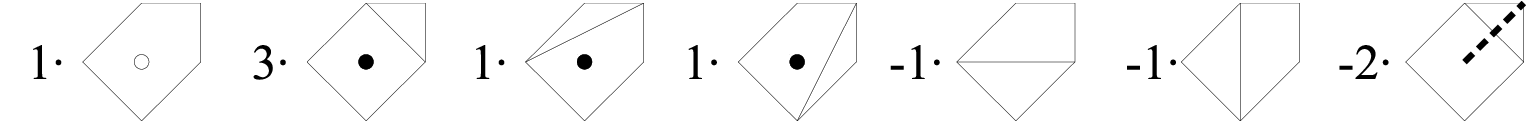}
\end{center}

Here the last picture denotes the cone that consists of all $\gamma$ such that $A^{\gamma}_0$ is the square and $A^{\gamma}_1$ is the dotted segment. This is the simplest example of a cone of a multisingularity stratum that cannot be represented as the set of all $\gamma\in\Z^A$ that induce a given regular subdivision of $A$. Although there is only one such cone in our example, most of cones are like this for large $A$. They are not among the cones of Severi varieties, whose weights are described in \cite{Kstud}.

The fan $(\pi_*\langle s=\frac{\partial s}{\partial x_{\cdot}}=0\rangle_{n+1})^2$ consists of the cones

\begin{center}
\noindent\includegraphics[height=1.2cm]{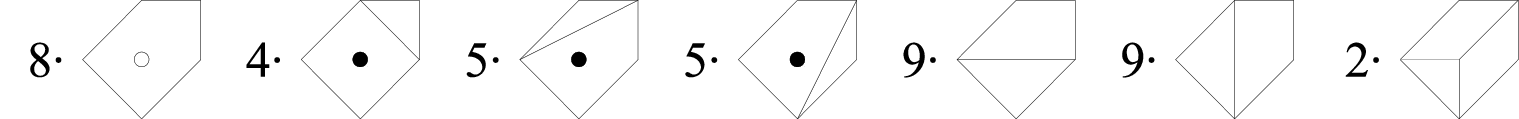}
\end{center}

See Remark \ref{remth3} on how to compute it. As a result, the tropical fan of the $\{\times\times\}$-stratum consists of the cones

\begin{center}
\noindent\includegraphics[height=1.2cm]{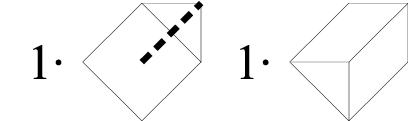}
\end{center}

and the tropical fan of the $\{\prec\}$-stratum consists of the cones

\begin{center}
\noindent\includegraphics[height=1.2cm]{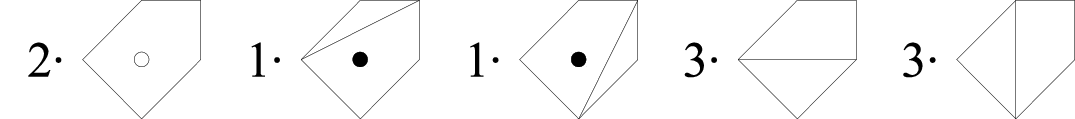}
\end{center}

which can be easily verified manually for this example.

{\noindent(A. Esterov) National Research University Higher School of Economics \newline Faculty of Mathematics NRU HSE, 7 Vavilova 117312 Moscow, Russia.}

\end{document}